\documentclass[times,final]{elsarticle}
\usepackage{jcomp}
\usepackage{framed,multirow}

\usepackage{amssymb}
\usepackage{latexsym}
\usepackage{url}
\usepackage[utf8]{inputenc} % allow utf-8 input
\usepackage[T1]{fontenc}    % use 8-bit T1 fonts
\usepackage{booktabs}       % professional-quality tables
\usepackage{amsfonts}       % blackboard math symbols
\usepackage{nicefrac}       % compact symbols for 1/2, etc.
\usepackage{microtype}      % microtypography
\usepackage{lipsum}
\usepackage{amsmath}
\usepackage{amsthm}
\usepackage{mathtools}
\usepackage{amsmath}
\usepackage{wrapfig}
\usepackage{comment}
\usepackage{todonotes}
\usepackage{caption}
\usepackage{xfrac}
\usepackage[]{algorithm2e}

\usepackage[arrow]{hhtensor}

\newtheorem*{remark}{Remark}
\newtheorem{theorem}{Theorem}
\newtheorem{example}{Example}

\newtheorem{proposition}[theorem]{Proposition}

\theoremstyle{definition}

\newcommand\restrict[1]{\raisebox{-.5ex}{$|$}_{#1}}

\newcommand{\jp}[1]{\ensuremath{[\![#1]\!]} }

\newcommand{\mhdbasis}{\vec{\mathcal{V}}}
\newcommand{\dgbasis}{\mathcal{V}}

\usepackage{mathtools}
\DeclarePairedDelimiter\abs{\lvert}{\rvert}%
\DeclarePairedDelimiter\norm{\lVert}{\rVert}%   

\makeatletter
\let\oldabs\abs
\def\abs{\@ifstar{\oldabs}{\oldabs*}}
\let\oldnorm\norm
\def\norm{\@ifstar{\oldnorm}{\oldnorm*}}
\makeatother

\journal{Journal of Computational Physics}

\begin{document}

\verso{Maria Han Veiga \textit{et al.}}

\begin{frontmatter}

\title{An arbitrary high-order Spectral Difference method for the induction equation}

\author[1]{Maria \snm{Han Veiga}\corref{cor1}}
\ead{mhanveig@umich.edu}
\cortext[cor1]{Corresponding author: }
%  Tel.: +0-000-000-0000;  
%  fax: +0-000-000-0000;}
\author[2]{David A. \snm{Velasco-Romero}}
\author[3]{Quentin \snm{Wenger}}
\author[2]{Romain \snm{Teyssier}}

\address[1]{Michigan Institute of Data Science,
University of Michigan, Weiser Hall
500 Church Street, Suite 600
Ann Arbor, MI 48109-1042}
\address[2]{Institute for Computational Science, University of Zurich, Winterthurerstrasse 190, 8057 Zurich, Switzerland}
\address[3]{ETH Zurich, Wolfgang Pauli Strasse 27,
8093 Zurich, Switzerland}

\received{xxx}
\finalform{xxx}
\accepted{xxx}
\availableonline{xxx}
\communicated{xxx}

\begin{abstract}
We study in this paper three variants of the high-order Discontinuous Galerkin (DG) method with Runge-Kutta (RK) time integration for the induction equation, analysing their ability to preserve the divergence-free constraint of the magnetic field. To quantify divergence errors, we use a norm based on both a surface term, measuring global divergence errors, and a volume term, measuring local divergence errors. This leads us to design a new, arbitrary high-order numerical scheme for the induction equation in multiple space dimensions, based on a modification of the Spectral Difference (SD) method \cite{Liu2006} with ADER time integration \cite{Dumbser2008}. It appears as a natural extension of the Constrained Transport (CT) method. We show that it preserves $\nabla\cdot\vec{B}=0$ exactly by construction, both in a local and a global sense. We compare our new method to the 3 RKDG variants and show that 
the magnetic energy evolution and the solution maps of our new SD-ADER scheme are qualitatively similar to the RKDG variant with divergence cleaning, but without the need for an additional equation and an extra variable to control the divergence errors. 
\end{abstract}

\begin{keyword}
%% MSC codes here, in the form: \MSC code \sep code
%% or \MSC[2008] code \sep code (2000 is the default)
%\MSC 41A05\sep 41A10\sep 65D05\sep 65D17
%% Keywords
%\KWD Keyword1\sep Keyword2\sep Keyword3
\end{keyword}

\end{frontmatter}

\section{Introduction}
\label{sec:introduction}

Developing numerical algorithms for the equations of ideal magneto-hydrodynamics (MHD) is of great interest in many fields of science, such as plasma physics, geophysics and astrophysics. Magnetic fields play an important role in a large variety of phenomena in nature, from the early universe, to interstellar and intergalactic medium, to environments and interiors of stars and planets \cite{Brandenburg2005}.

The ideal MHD equations describe conservation laws for mass, momentum and total energy on the one hand, and for magnetic flux on the other hand. The first 3 conservation laws form what is called the Euler sub-system, while the fourth one is called the induction sub-system. In this paper, we focus  on the latter, usually called the induction equation:

\begin{equation}
\label{eq:mhd-induction-eq}
\partial_t \vec{B} =  \nabla\times(\vec{v}\times \vec{B}) + \eta \nabla^2 \vec{B}.\end{equation}

This partial differential equation describes the evolution of a magnetic field $\vec{B}$ under the effect of the velocity $\vec{v}$ of an electrically conductive fluid. The coefficient $\eta$ denotes the magnetic diffusivity. In the ideal MHD case, the fluid has infinite electric conductivity, so that $\eta \to 0$ and the diffusive term can be ignored. 

By taking the divergence of Eq.~\eqref{eq:mhd-induction-eq}, we note that the time evolution of the divergence of $\vec{B}$ is zero for all times, meaning that the initial divergence of $\vec{B}$ is preserved:

\begin{equation}
\label{eq:divbevol}
\partial_t (\nabla \cdot \vec{B})  = \nabla \cdot \left( \nabla\times(\vec{v}\times \vec{B}) \right) = 0,
\end{equation}
as the divergence of the curl of a vector is always zero.

Physically, the fact that magnetic fields have no monopoles and that magnetic field lines form closed loops, is translated in the initial condition
\begin{equation}
\label{eq:divfree}
\nabla \cdot \vec{B} = 0.
\end{equation}
Considering Eq.~\eqref{eq:divbevol} and Eq.~\eqref{eq:divfree} together means that the divergence of $\vec{B}$ must be zero at all times.
To clearly see the erroneous evolution of our system if $\nabla\cdot\vec{B}$ happens to be nonzero, we can re-formulate Eq.~\eqref{eq:mhd-induction-eq} as
\begin{equation}
\label{eq:induction-eq-ext}
\partial_t \vec{B} + (\vec{v}\cdot\nabla)\vec{B} = -\vec{B}(\nabla\cdot \vec{v}) + (\vec{B}\cdot\nabla)\vec{v} + \vec{v}(\nabla\cdot \vec{B}).
\end{equation}
Note that the second term of the left-hand side, $(\vec{v}\cdot\nabla)\vec{B}$, corresponds to the advection of $\vec{B}$ by the fluid, the first term on the right-hand side models the \textit{compression} of the magnetic field lines and the second term is due to the \textit{stretching} of the field lines. This interpretation can be done by establishing an analogy with the vorticity equation \citep[see][for example]{davidson2001}. The last term, proportional to $\nabla\cdot \vec{B}$, is also proportional to the velocity of the flow $\vec{v}$, and vanishes only if the magnetic field is divergence free.

When applying common discretisation schemes, for example the popular Finite Volume (FV) method, the divergence-free constraint in Eq.~\eqref{eq:divfree} is not necessarily fulfilled in the discrete sense. Indeed, in this case, the numerical representation of the field is based on a volume integral of the magnetic field and the magnetic flux is not conserved anymore. 

This is a big issue in the numerical evolution of the MHD equations, as shown for example in the seminal studies of \cite{BRACKBILL1980,toth1996}, which show that a non-physical force parallel to the velocity $\vec{v}$ and proportional to $\nabla \cdot \vec{B}$ appears in the discretised conservative form of the momentum equation in the Euler sub-system.

There have been many proposed methods to guarantee a divergence-free description of $\vec{B}$. For example, the non-solenoidal component of $\vec{B}$ is removed through a Hodge-Helmholtz projection at every time step (e.g. \cite{BRACKBILL1980,zachary1994}), or the system in Eq.~\eqref{eq:induction-eq-ext} is written in a non-conservative formulation where the non-solenoidal component of $\vec{B}$ is damped and advected by the flow (e.g. \cite{powell1999,dedner2002,munz1999}). 

Another approach is done at the discretisation level, where the numerical approximation of the magnetic field is defined as a surface integral and collocated at face centres, while the electric field used to updated the magnetic field is collocated at edge centres, in a staggered fashion \cite{yee1966, brecht1981, Evans1988, devore1989}. This method, called Constrained Transport (CT), was later adapted to the FV framework applied to the MHD equations \citep[see e.g.][]{Dai1998, Ryu_1998, Balsara_mhd_1999, Balsara_2004, Fromang2006}. The CT method is obviously closer to the original conservation law, as the magnetic flux is explicitly conserved through the cell faces. A comprehensive review of these methods in astrophysics can be found in \cite{Teyssier2019} and references therein.

In addition, finite element methods can be naturally used to solve the induction and MHD type equations. In particular, when the magnetic field is approximated by a $H(div)$ vector function space (where elements of this space have square integrable divergence), it leads to continuous normal components of the approximation across element faces, while when the electric field is approximated by a $H(curl)$ vector function space (elements of this space have square integrable curl), it leads to continuous tangential components across cell faces \cite{Brezzi_1991}. For example, Raviart-Thomas/N\'ed\'elec basis functions are conforming with the aforementioned vector function spaces and have been used successfully to solve the induction equation \cite{balsara_kappeli_2018, chandrashekar_2020, praveen_2019}. %Traditionally, because there is a requirement on the continuity across the cell faces, global matrices are used. However, \cite{praveen_2019} presents a DG-like method based on Raviart-Thomas elements which remains local.}

With the increased availability of high-order methods, one could ask whether a high-order approximation of the magnetic field $\vec{B}$ alone could be sufficient to control the non-vanishing divergence problem of the magnetic field. Indeed, very high-order methods have been developed in many fields of science and engineering, both in the context of the FV method \cite{Jiang1999} and in the context of the Finite Element (FE) method \cite{Li2005,Mocz2013,Fu2018,Guillet2019}. These very high-order methods have proven successful in minimizing advection errors in case of very long time integration \cite{Gassner2013,Sengupta2006,Velasco2018}. Very high-order methods have already been developed specifically for the ideal MHD equations 
\citep{Nordlund1990,Jiang1999,Balsara_2004,Balsara_weno_2009,Felker2018,balsara_kappeli_2018}. It turns out, as we also show in this paper, that a very high-order scheme does not solve by itself the problem of non-zero divergence and specific schemes have to be developed to control the associated spurious effects \citep{munz1999,Li2012,Fu2018,Guillet2019}. 

In this paper, we present a new, arbitrary high-order method that can perform simulations of the induction equation, based on the Spectral Difference method developed in \cite{Kopriva1998, Liu2006} and on the ADER timestepping scheme 
\cite{dumbser_ader_2013, balsara_ader_2018}. We show that this technique is by construction strictly divergence free, both in a local and in a global sense. While there are similarities between this work and the work presented in \cite{balsara_kappeli_2018}, there are some key differences: our scheme includes internal nodal values which are evolved according to a standard SD scheme, similar to \cite{praveen_2019}. Furthermore, there is no need for an explicit divergence-free reconstruction step, which means achieving arbitrarily high-order is simpler. In particular, Propositions \ref{proposition:pointwise_div_free} and \ref{proposition:globally_div_free} (see below) prove that our method is divergence-free by construction and arbitrarily high-order.

The paper is organised as follows: we start in section \ref{sec:overview} with a detailed description of several well-known high-order methods used to model the induction equation, discussing the challenges of controlling efficiently the magnitude of the divergence of the magnetic field. Then, in section \ref{sec:sd_mhd}, we present our new Spectral Difference method for the induction equation, highlighting our new solution points for the magnetic field components and the need of two-dimensional Riemann solvers. In section \ref{sec:numerics} we evaluate the performance of the new SD-ADER method numerical through different test cases. In section \ref{sec:mhd-discussion}, we compare our new method to other very high-order schemes using different numerical experiments. Finally,  in section \ref{sec:conclusion}, we present our conclusions and outlook.
\section{Various high-order Discontinuous Galerkin methods for the induction equation}
\label{sec:overview}

\begin{figure}
   \includegraphics[width=0.94\textwidth]{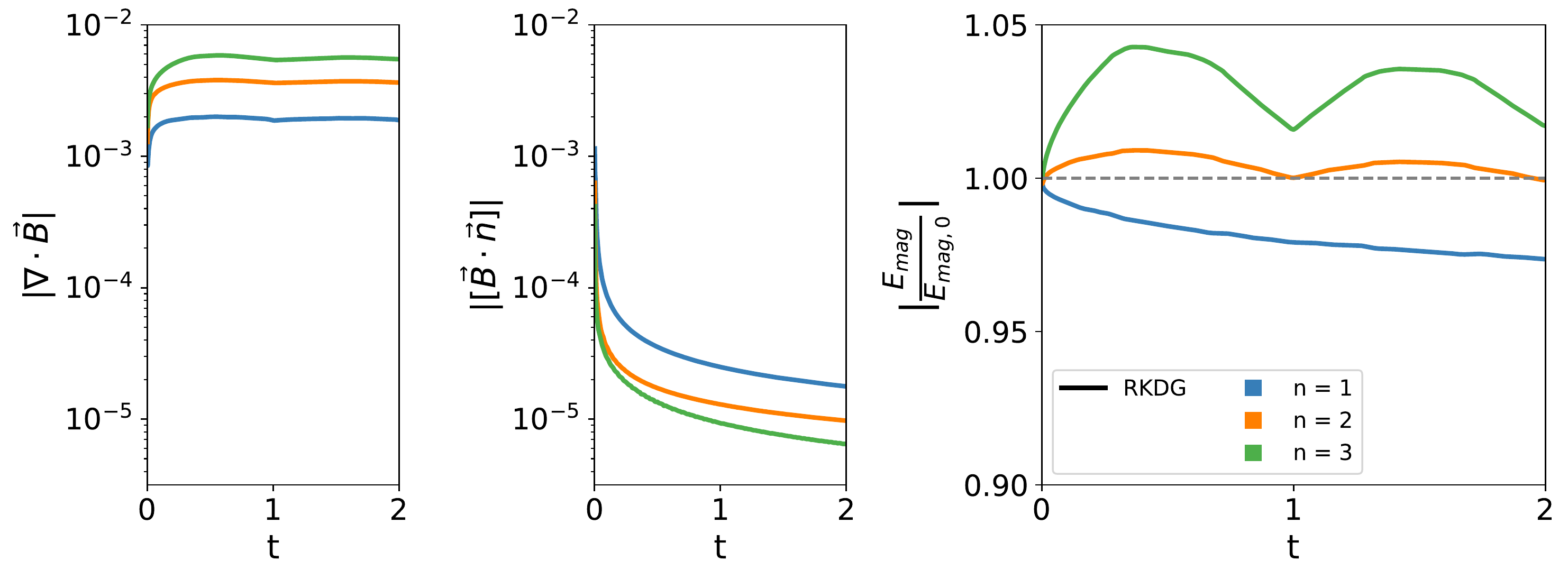}
   \begin{center}
   \includegraphics[width=1.0\textwidth]{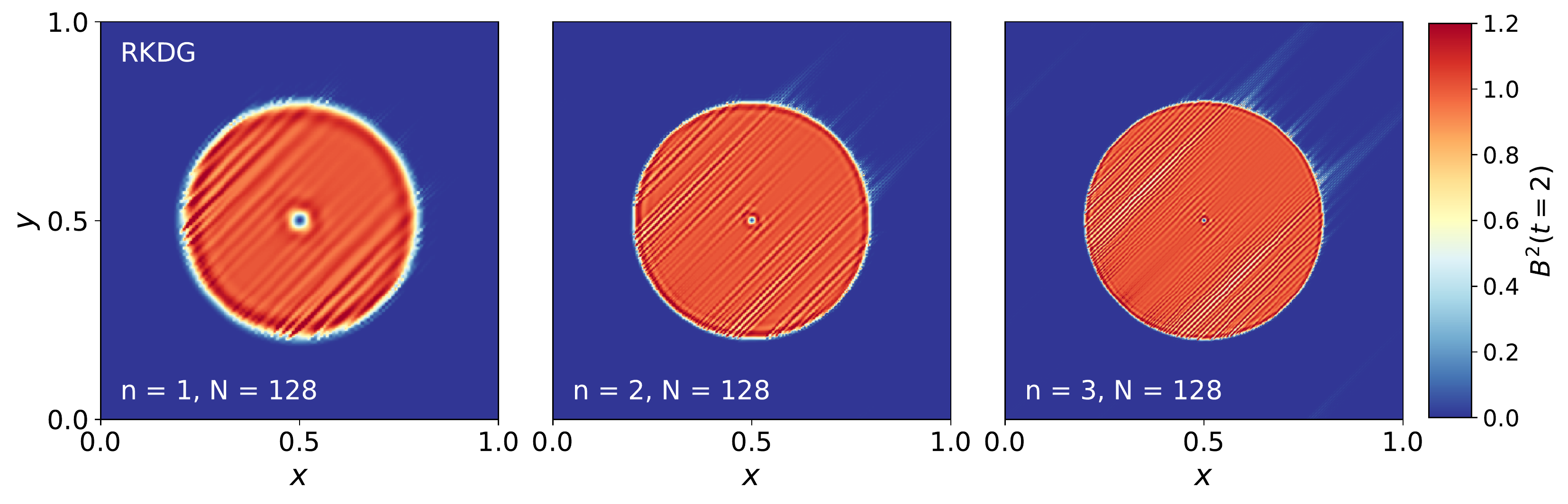}
    \caption{Performance of the traditional RKDG scheme for the magnetic loop test defined in Eq.~\eqref{eq:magloop-ics}. In the top row, the first two panels show the divergence contribution of the volume term and of the surface term, respectively. The third panel shows the magnetic energy of the solution over time. In the bottom row, maps of the magnetic energy density are shown at $t = 2$. The three runs corresponds to increasing polynomial degree ($n=1$, $n=2$ and $n=3$) and a fixed number of cells ($N = 128$ per dimension).  \label{fig:rkdg-loopadvection-div-energy}}
    \end{center}
\end{figure}

The Discontinuous Galerkin (DG) method is a very popular FE scheme in the context of fluid dynamics. It is based on a weak formulation of the conservative form of the MHD equations using a Legendre polynomial basis \cite{cockburn1998}. In this context, the induction equation has to be written as 
\begin{equation}
\partial_t \vec{B} +  \nabla \cdot (\vec{v}\otimes \vec{B} - \vec{B}\otimes \vec{v}) = 0,
\label{eq:mhd-induction-eq-div}
\end{equation}
in a conservative form compatible with the Euler sub-system. Indeed, this equation is now based on the divergence operator, which forces us to deal with the magnetic field through volume integrals.

In this section, we describe three different implementations of the DG method for the induction equation. The first one is the classical DG method with Runge-Kutta time integration (RKDG), for which nothing particular is done to preserve the divergence-free character of the solution. The second method, presented for example in \cite{Guillet2019,Cockburn2004}, uses a modified polynomial basis for the magnetic field, so that it is locally exactly divergence free. The third one allows the divergence to explicitly deviate from zero, but tries to damp the divergence errors using an additional scalar field and its corresponding equation \citep{munz2000}. We will evaluate the performance of these three classical methods using a proper measurement of the divergence error, as well as the conservation of the magnetic energy, using the famous magnetic loop advection test. 

\subsection{A traditional RKDG for the induction equation}

In this section we describe the classical modal RKDG method using a simple scalar problem in one space dimension: 
\begin{equation}
\label{eq:conslaw}
\begin{cases}
      \partial_t u + \partial_x f( u ) = 0 \quad \in \Omega \times [0,\infty]\\
      u(t=0) = u_0 \\
      u_{\partial \Omega} = g.
    \end{cases}
\end{equation}
The generalisation to multiple space dimensions for structured Cartesian grids can be achieved through tensor products.
Let $\Omega \in \mathbb{R}$ be a regular domain which is discretised by $N$ elements $K_p =  [x_{p-1/2},x_{p+1/2}]$ for $p=1,...,N$.
Consider the local space $\dgbasis$ given by the set $\{\phi_i\}_{i=0}^{n}$ of one dimensional Legendre polynomials with degree of at most $n$ in $x$. 
For each element $K_{p}$, the numerical solution is written as:
\[u(x,t) = \sum_{i=0}^{n} \hat{u}_i(t) \phi_i(x),\]
where the modal coefficient $\hat{u}_i(t)$ is obtained by the $L^2$ projection of the solution $u(x)$ on the $i$-th Legendre basis polynomial.
The DG method is based on a weak form of Eq.~\eqref{eq:conslaw}, projecting it on the polynomial basis, followed by an integration by parts.
We obtain the following semi-discrete formulation of the DG method as:
\begin{align*}
\label{eq:dg}
\frac{d \hat{u}_i}{dt} + \left[ \hat{f}(u(x,t))\phi_i(x)\right]_{x_{p-1/2}}^{x_{p+1/2}} - \int_{K_p} f(u(x,t)) \partial_x \phi_i(x) {\rm d}x = 0,\quad i=0,...,n,
\end{align*}
where we exploited the fact that Legendre polynomials form an orthonormal basis. Note that the surface term in the previous equation needs a Riemann solver to compute a continuous
numerical flux at element boundaries, noted here $\hat{f}$.
Once the spatial component has been discretised, we are left with an ordinary differential equation of the form:
\[ \frac{d}{dt} u = \mathcal{L}(u), \]
where $\mathcal{L}$ denotes the DG discretisation operator. Integration in time is performed using a Strong Stability Preserving (SSP) RK method \cite{Gottlieb2005, Kubatko2014}. The time step has to fulfill a Courant-Friedrich-Lewy (CFL) condition to achieve numerical stability, which for the RKDG scheme reads \cite{cockburn1998}:
\[
\Delta t = \frac{C}{2n + 1} \frac{\Delta x}{\left| v_{\rm max}\right| },
\]
where $n$ is the polynomial degree and $C$ is a constant usually set to $C=0.8$.

\subsection{Quantifying divergence errors}

It is highly non-trivial to estimate the error in the divergence of the magnetic field for high-order methods in general, and for FE schemes in particular. Indeed, the numerical approximation of the solution is defined in a local sense, with polynomials of degree at most $n$ inside each element, but also in a global sense by considering the solution given by the union of all the elements.
A suitable measurement for $\nabla\cdot \vec{B}$ has been proposed by \cite{Cockburn2004} as
\begin{equation}
\label{eq:divmeas}
\norm{\nabla \cdot \vec{B}} = \sum_{e\in \mathcal{E}} \int_e \abs{ \jp{ \vec{B}\cdot \vec{n} } } {\rm d}s + \sum_{K\in \mathcal{K}} \int_K \abs{\nabla \cdot \vec{B}} {\rm d}\vec{x},
\end{equation}
where $\jp{ \vec{B}\cdot\vec{n}_x } = B_x^{int(K)}-B_x^{ext(K)} $ (for example) denotes the jump operator and $B_x^{int(K)},~ B_x^{ext(K)}$, are the limits of $B_x$ at interface $e$ from the interior and exterior of $K$ respectively. We assume $\vec{B}$ is smooth within each element $K \in \Omega$.
However, in the DG framework, $\vec{B}$ can be discontinuous across element boundaries (noted here $e$). In the previous formula, $\mathcal{E}$ denotes the set of element interfaces and $\mathcal{K}$ the set of element volumes. Note that, for a piecewise-smooth function that is divergence free inside each element, it is globally divergence free if and only if the normal component of the vector field across each interface $e$ is continuous, hence the consideration of the jump in the normal component of the magnetic field across the interfaces $e$, given by the first term in Eq.~\eqref{eq:divmeas}. This divergence error measurement has been derived by exploiting the properties of the space $H(div)$ \cite{nedelec1980} or by using a functional approach \cite{Cockburn2004}. In what follows, we call the first contribution the surface term, and the second contribution the volume term.

\subsection{Magnetic energy conservation}

The other metric used in this paper to characterise different numerical methods is the evolution of the magnetic energy. This is particularly important in the context of magnetic dynamos \citep{Roberts2000, Brandenburg2005}, as one wishes to avoid having spurious magnetic dynamos triggered by numerical errors. Using \eqref{eq:induction-eq-ext} and considering again a non-zero divergence, the magnetic energy equation can be written as:
\begin{equation}
\label{eq:induction-energy-eq}
\partial_t \left( \frac{B^2}{2}  
\right) 
+ \left( \vec{v}\cdot \nabla \right) \left( \frac{B^2}{2} 
\right) 
= - B^2(\nabla\cdot\vec{v}) + \vec{B}\cdot(\vec{B}\cdot\nabla)\vec{v} + (\vec{B}\cdot\vec{v})(\nabla\cdot\vec{B}),
\end{equation}
where the last term is here again spurious.

For example, in the simple case of pure advection where $\vec{v}$ is constant, one can observe that the first two terms on the right hand side vanish, while the third term vanishes only if $\nabla\cdot\vec{B}=0$. 
On the other hand, if $\nabla\cdot\vec{B} \ne 0$, depending on the solution properties, one could observe a spurious increase of the magnetic energy over time, and interpret it wrongly as a dynamo. In the advection case, the magnetic energy is expected to remain constant, although we expect the numerical solution of the magnetic energy to decay, owing to the numerical dissipation associated to the numerical method. It should however never increase.
    
\subsection{The field loop advection test}    

The advection of a magnetic loop is a well-known numerical experiment introduced for example in \cite{Gardiner2005} to test
the quality of the numerical solution, with respect to both divergence errors and magnetic energy conservation.
The test is defined using the following {\it discontinuous} initial magnetic field,
\begin{equation}
\label{eq:magloop-ics}
    \vec{B}_0 = \begin{pmatrix} B_{x,0} \\ B_{y,0} \end{pmatrix} =  \begin{pmatrix} -A_0(y-y_c)/r \\ A_0(x-x_c)/r \end{pmatrix} \quad {\rm ~for~}r < r_0,
\end{equation} 
and $\vec{B}_0=0$ otherwise, advected with a constant velocity field $\vec{v}=(1,1)$. We use here $A_0 = 0.001$, $r_0=0.25$ and $(x_c,y_c)=(0.5, 0.5)$. We consider a square box $[0,1]\times[0,1]$ and the final time $t = 2$. This allows the loop to cross the box twice before returning to its initial position.

In Fig.~\ref{fig:rkdg-loopadvection-div-energy} we show the performance of our traditional RKDG scheme at different approximation orders. When measuring the divergence errors of the numerical solution, we observe that the volume term (measuring global divergence errors) seems to decrease with the approximation order (middle panel) as expected. On the contrary, the surface term (measuring local divergence errors) does not decrease at all. In fact, local errors increase with increasing polynomial degree. 

Furthermore, the magnetic energy evolution is clearly incorrect. Namely, at $3^{rd}$ and $4^{th}$ orders (corresponding to a maximal polynomial degree of $n=2$ and $n=3$, respectively), an initial increase on the magnetic energy is observed. In the bottom panel (Fig.~\ref{fig:rkdg-loopadvection-div-energy}), maps of the magnetic energy density $B^2/2$ (normalised to the maximum value in the initial condition) are shown at $t = 2$ and at different orders. We see spurious stripes with high-frequency oscillations, aligned with the direction of the velocity field. Our results are similar to the numerical experiments performed in  \cite{nunez2018} and consistent with Eq.~\eqref{eq:induction-energy-eq}. We clearly have to give up our initial hopes that going to very high order would solve the divergence-free problem.

\subsection{RKDG with a locally divergence-free basis (LDF)}
\label{sec:ldfrkdg}
The locally divergence-free (LDF) method was first introduced by \cite{Cockburn2004} with the intention to control the local contribution of the divergence. Indeed, we have seen in the last sub-section that this term dominates the error budget in our simple numerical experiment. This method has been recently revisited in \cite{Guillet2019} in conjunction with several divergence cleaning schemes.

LDF is built upon the previous RKDG scheme with the key difference that the approximation space considered for the magnetic field $\vec{B}$ is given by:
\[ \mhdbasis^n = \{ \vec{v} \in [L^1]^d: \vec{v}\restrict{K} \in[\mathbb{P}^n]^d, \nabla \cdot \vec{v}\restrict{K} = 0 \}. \]
The trial space considered contains only functions which are divergence free inside each element $K$ and belong to the $d$-dimensional vector space  $[\mathbb{P}^n]^d$, where each polynomial is a polynomial of at most degree $n$. One key difference between this method and the traditional RKDG is that the modal coefficients of the solution are now shared between $B_x, B_y$ and $B_z$ due to this new carefully designed vector basis. We show in this paragraph only the example for $n=1$ in two space dimensions. For more details on the implementation, please refer to Appendix~\ref{ap:LDF}.

\begin{example}
{\textbf{d = 2, n = 1:}}
Consider the basis elements of the $\mbox{span}(\{1,x,y,xy,y^2,x^2\}) = \mathbb{P}^{2}(x,y)$. Form the vector $\vec{b}_i = (0, 0, v_i)^T$ for $v_i \in {\rm basis}(\mathbb{P}^{2}(x,y))$ and take the $curl$ of its elements. This set of vectors spans a subspace of $[\mathbb{P}^{1}(x,y)]^2$.
\[ \mhdbasis^1 = {\rm span}\left(\{(0,-1),(1,0),(-x,y),(2y,0),(0,2x)\}\right) \subset [\mathbb{P}^1(x,y)]^2. \]
\end{example}

\begin{figure}
   \includegraphics[width=0.94\textwidth]{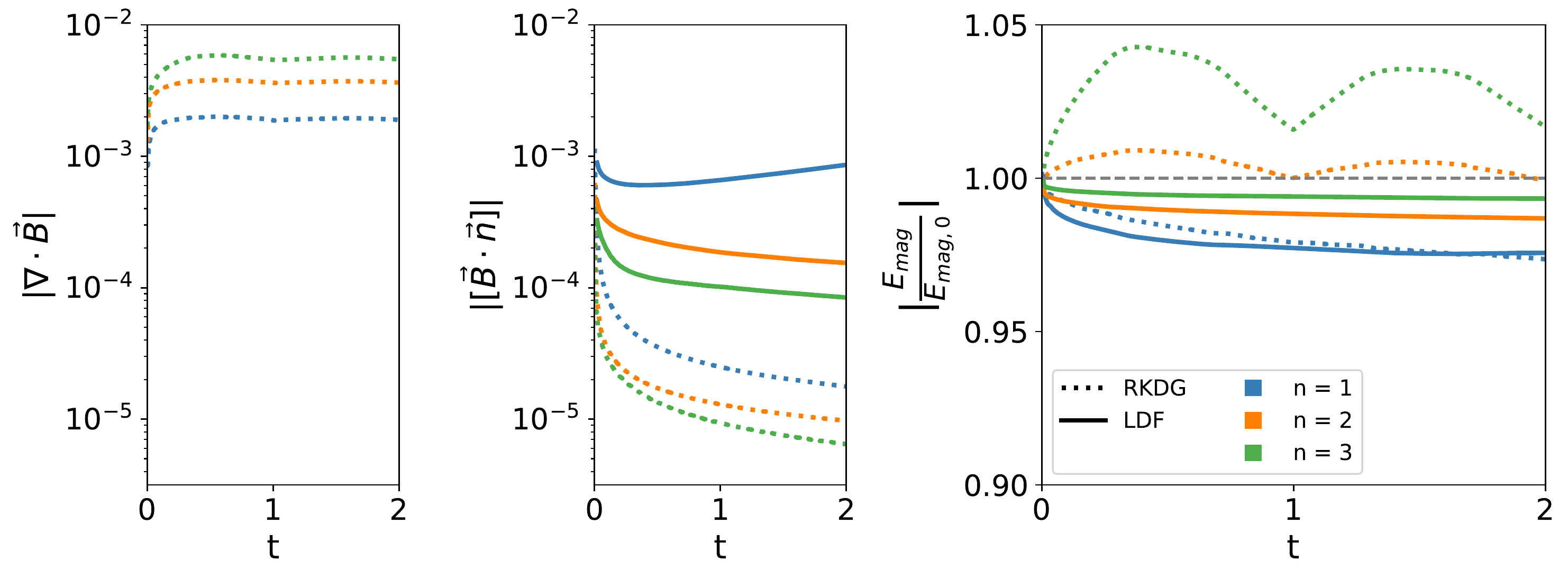}
   \begin{center}
   \includegraphics[width=1.0\textwidth]{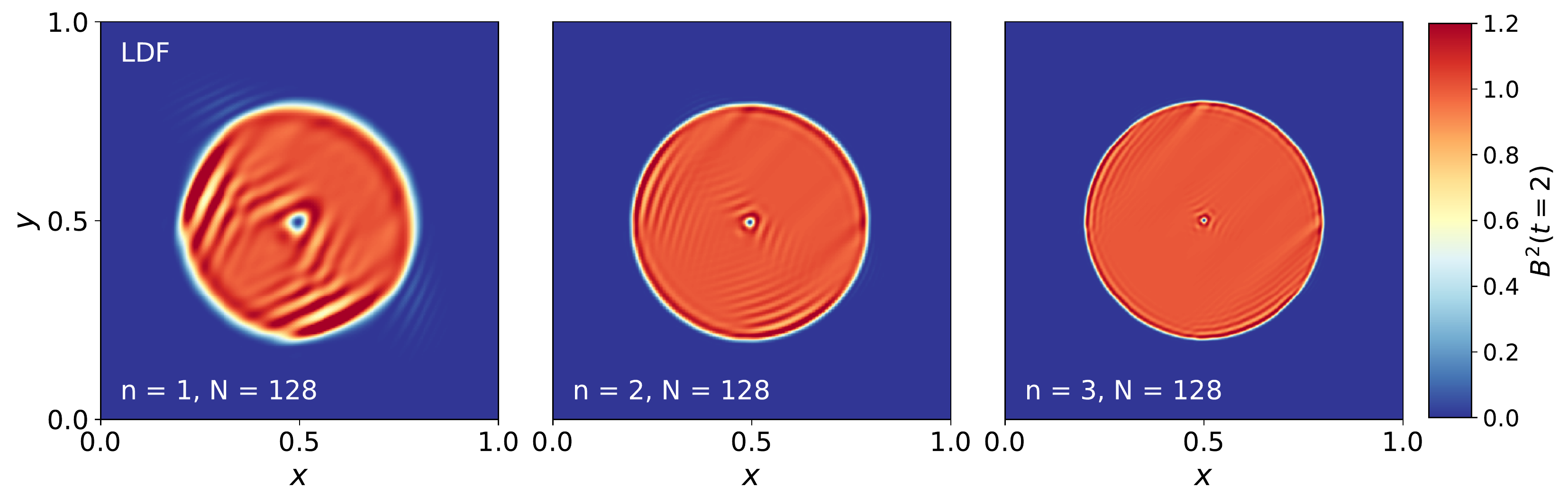}
    \caption{
    Same as Fig.~\ref{fig:rkdg-loopadvection-div-energy} but now for the LDF scheme (solid lines). For comparison, the results of the traditional RKDG scheme are shown as dotted lines. Note that the LDF scheme has no local divergence errors by construction (no solid lines in the top left panel). \label{fig:ldf-loopadvection-div-energy}}
    \end{center}
\end{figure}

In Fig.~\ref{fig:ldf-loopadvection-div-energy} we show the performance of the LDF scheme at different approximation orders. When measuring the divergence of the numerical solution, we observe that the local contribution of the divergence is zero (as expected). The global contribution (middle panel), while decreasing with the order, is considerably larger than the traditional RKDG scheme. We believe this is due to the reduced number of degrees of freedom in the LDF basis. For the measured magnetic energy, we don't see a spurious dynamo anymore, only a decay due to numerical dissipation. We also observe less numerical dissipation when increasing the order of the method, a desirable property. In the bottom panel of the same figure, the magnetic energy density maps show some residual oscillations at $t = 2$, although much less than for the original RKDG scheme. In order to reduce these oscillations even more, one traditionally uses the LDF scheme in conjunction with a divergence cleaning strategy \cite{Cockburn2004,Guillet2019,klingenberg2017}, similar to the one presented in the next section.

\subsection{RKDG with hyperbolic and parabolic divergence cleaning (DivClean)}
Divergence cleaning is a general strategy aiming at actively reducing divergence errors by modifying the solution at every time step. Among many possibilities that can be found in the literature \citep[see][for example]{toth1996}, we have adopted here a robust technique based on the addition of a new variable that can be used to control and dissipate the divergence of the magnetic field $\vec{B}$. Following \cite{dedner2002}, we briefly describe this method that performs what is called parabolic and hyperbolic \textit{divergence cleaning}. The idea is to introduce an additional scalar field $\psi$ and couple it to the induction equation. This method is also known as the Generalised Lagrangian Multiplier (GLM) approach \cite{munz1999, munz2000}.

The induction equation in its divergence form in Eq.~\eqref{eq:mhd-induction-eq-div} is modified as
\begin{equation}
\label{eq:glm-induction-eq}
\begin{split}
\partial_t \vec{B} + \nabla \cdot (\vec{v}\otimes \vec{B} - \vec{B} \otimes \vec{v}) + \nabla \psi &= 0,\\
\mathcal{D}(\psi) + \nabla \cdot \vec{B} &= 0,
\end{split}
\end{equation}
where $\mathcal{D}(\cdot)$ is a linear differential operator. There are different ways to choose $\mathcal{D}(\cdot)$ \cite{munz1999,munz2000}. In this work, we choose a \textit{mixed} type of correction, defining  $\mathcal{D}(\cdot)$ as

\[ \mathcal{D}(\psi):= \frac{1}{c_h^2}\partial_t \psi + \frac{1}{c_p^2}\psi. \]

The new scalar variable $\psi$ is coupled to the non-vanishing divergence of the magnetic field and evolves according to a new additional partial differential equation:

\[\partial_t \psi + \frac{c_h^2}{c_p^2}\psi + c_h^2\nabla\cdot\vec{B} = 0.\]

Both $c_h$ and $c_p$ are free parameters tuned for each particular problem at hand. 
The hyperbolic parameter $c_h$ corresponds to the velocity of the waves that are carrying the divergence away from regions where errors are created. The parabolic parameter $c_p^2$ corresponds to the diffusion coefficient of the parabolic diffusion operator that damps divergence errors.
There are different strategies to choose $c_h$ and $c_p$ that could lead to a robust scheme. Different methods have been proposed in the literature \cite{dedner2002,Guillet2019,Mignone2010}, and these choices boil down to setting the speed $c_h$ to be a small multiple of the maximum of the velocity field $\left| v_{\rm max} \right|$ and the magnitude of the diffusion coefficient $c_p^2$ as a small multiple of $c_h \Delta x$.

\begin{figure}
   \includegraphics[width=0.94\textwidth]{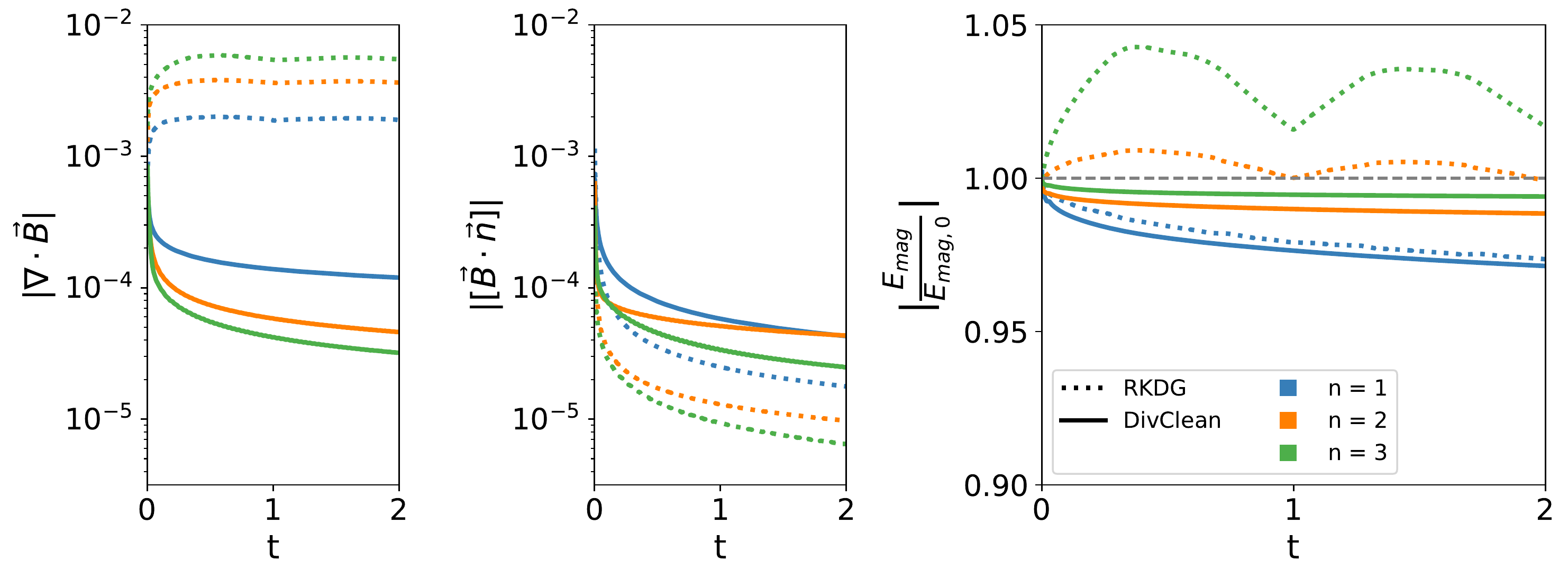}
   \begin{center}
   \includegraphics[width=1.0\textwidth]{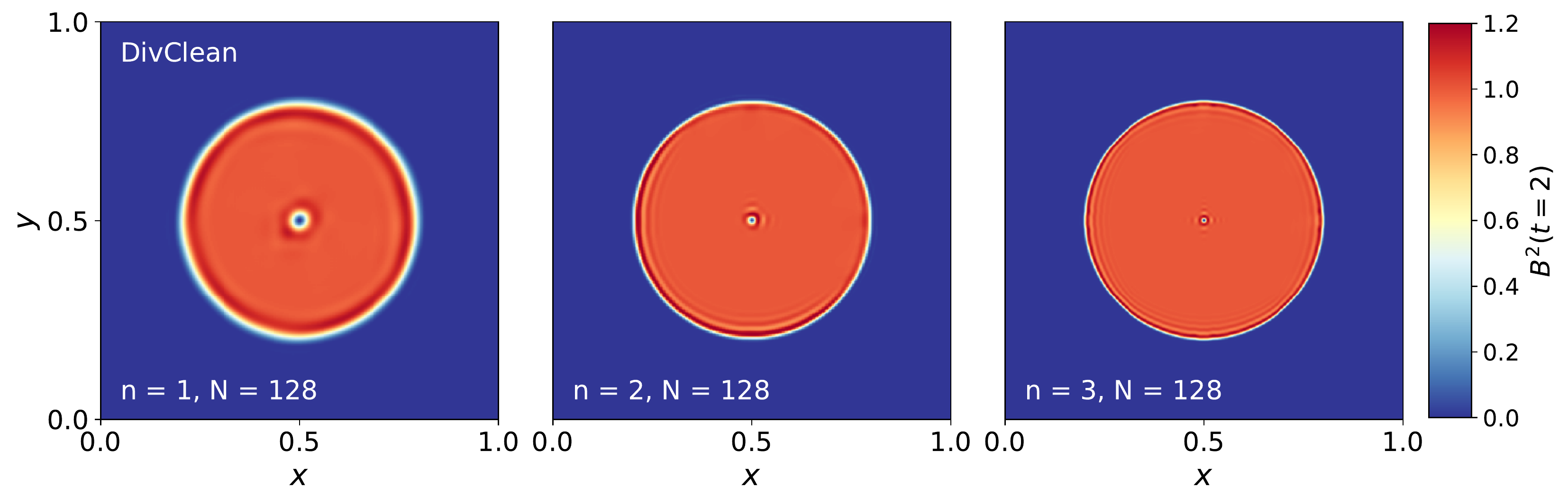}
    \caption{Same as Fig.~\ref{fig:rkdg-loopadvection-div-energy} but now for the DivClean scheme (solid lines). for comparison, the results of the traditional RKDG scheme are shown as dotted lines. \label{fig:divc-loopadvection-div-energy}}
    \end{center}
\end{figure}

In Fig.~\ref{fig:divc-loopadvection-div-energy} we show the performance of the RKDG scheme with both hyperbolic and parabolic divergence cleaning, called here DivClean, at different approximation orders. For implementation details, please refer to Appendix~\ref{ap:DivClean}.  In this numerical experiment, we set $c_h$ and $c_p^2$ according to \cite{Mignone2010}, namely, we choose $c_h = 2 \left| v_{\rm max}\right| $ such that the overall time step is almost not affected, and $c_p^2 = 0.8 c_h \Delta x$. We see that both surface and volume terms of the divergence error norm are small, and they both decrease with increasing orders. The magnetic energy density maps look very smooth and symmetrical, with very small residual features close to the discontinuity. It is worth stressing that none of the tests performed here make use of TVD slope limiters, so that some residual oscillations due to the Runge phenomenon are expected.

\section{A high-order Spectral Difference method with Constrained Transport}
\label{sec:sd_mhd}

We present in this section a new method within the FE framework that addresses most of the issues we discussed in the previous section. It is both locally and globally divergence free, and it does not need the introduction of a new variable and a new equation, together with free parameters sometimes difficult to adjust. This new method is based on the Spectral Difference (SD) method \cite{Liu2006}. In this section, we present the original SD method for the semi-discrete scheme, followed by a description of
our time integration strategy. We then focus on the modifications to the original method to solve the induction equation.
We prove that our method is both strictly divergence free at the continuous level inside the elements, and at the interface between elements,
maintaining the strict continuity of the normal component of the field. Using Fourier analysis, we finally show  that our new method attains the same stability properties as reported in \citep{abeele2008,jameson2010}. %is stable for standard control volume discretisation strategies.

\subsection{The SD method in a nutshell}

For sake of simplicity, we present the SD method using a simple scalar problem in one space dimension. The generalisation to multiple
space dimensions will be discussed later.
Let us denote the numerical solution as $u(x)$. 
We focus on the description of the solution in one element, which is given by Lagrange interpolation polynomials $\{\ell^s_{i}(x)\}_{i=0}^n$, built on a set of points $\mathcal{S}^s=\{x^s_i\}_{i=0}^n$, called the solution points (with the superscript $s$). The numerical solution inside an element is given by:
\[ u(x) = \sum_{i=0}^n u(x^s_i)\ell^s_i(x), \]
where $n$ is the polynomial degree of the interpolation Lagrange polynomials.
The SD method features a second set of nodes $\mathcal{S}^f=\{x^f_i\}_{i=0}^{n+1}$, called the flux points (with the superscript $f$). A numerical approximation of the flux is evaluated by another set of Lagrange interpolation polynomials $\{\ell^f_{i}(x)\}_{i=0}^{n+1}$ built on the flux points. Note that we have $n+1$ solution points and $n+2$ flux points, and that the first and the last flux points coincide with the boundary of the elements ($x^f_0$ and $x^f_{n+1}$). Moreover, at the interfaces between elements, a numerical flux based on a Riemann solver must be used to enforce the continuity of the flux between elements. Let $\hat{f}(\cdot)$ denote this single-valued numerical flux, common to the element and its direct neighbour. The  approximation for the flux is given by:
\begin{equation}
\label{eq:numerical_flux}
f(x) = \hat{f}(u(x^f_0)) \ell^f_0(x) + \sum_{i=1}^{n} f_i(u(x^f_i)) \ell^f_i(x) +  \hat{f}(u(x^f_{n+1})) \ell^f_{n+1}(x),
\end{equation}
where we wrote separately the two extreme flux points with their corresponding  numerical flux.
The final update of the solution is obtained using the exact derivative of the flux evaluated at the solution points, so that the semi-discrete scheme reads:
\begin{align*} 
\frac{{\rm d}}{{\rm d}t} u(x^s_j) = - \hat{f}(u(x^f_0))  \ell^{f\prime}_0(x^s_j) -\sum_{i=1}^{n} f_i(u(x^f_i)) \ell^{f\prime}_i(x^s_j) - \hat{f}(u(x^f_{n+1})) \ell^{f\prime}_{n+1}(x^s_j),
\end{align*}
where the primes stand for the derivative of the Lagrange polynomials. A straightforward extension to more space dimensions can be achieved by the tensor product between the set of one dimensional solution points and flux points. The left panel on Fig.~\ref{fig:sd_representation} shows in blue the solution points and in red (and salmon) colour the flux points for a classical SD scheme in two space dimensions, as well as the subcells (denoted by the black lines) which we call \textit{control volumes}. The stability of the SD method in one dimension has been shown in \cite{jameson2010} at all orders of accuracy, while the stability of the SD scheme in two dimensions, for both Cartesian meshes and unstructured meshes, has been demonstrated in \cite{abeele2008}. 

As shown in \cite{jameson2010}, the stability of the standard SD method depends on the proper choice of the flux points and not on the position of the solution points. The only important requirement is that the solution points must be contained anywhere within the (inner) control volume delimited by the flux points. With this in mind, we use Gauss-Legendre quadrature points for the inner flux points and the zeros of the Chebyshev polynomials for the solution points, and we show in section \ref{sec:stability} that indeed this general result also holds for the induction equation.

\subsection{High-order time integration using ADER}

We decided not to use the same SSP Runge Kutta method as for the RKDG scheme. 
Instead, we decided to explore the modern version of the ADER method \cite{Dumbser2008,balsara2009,mhv2020}.
Indeed, we believe this method is well suited to compute solutions to arbitrary high order in time.
We exploit this nice property in our numerical experiments shown in section \ref{sec:numerics}.
Consider again the scalar, one-dimensional conservation law given in Eq.~\eqref{eq:conslaw},
\begin{equation}
\begin{cases}
      \partial_t u + \partial_x f( u ) = 0 \quad \in \Omega \times [0,\infty]\\
      u(t=0) = u_0 \\
      u_{\partial \Omega} = g,
    \end{cases}
\end{equation}
with suitable initial conditions and boundary conditions. For simplicity, we are only updating the solution $u(x^s_i,t)$ for a single solution point $x^s_i$. 
Modern ADER schemes are based on a Galerkin projection in time.  We multiply the previous conservation law by an arbitrary test function $\psi(t)$, 
integrating in time over $\Delta t$:
\[\int^{\Delta t}_0  \psi(t)\partial_t u {\rm d}t + \int^{\Delta t}_0 \psi(t)\partial_x f(u)  {\rm d}t = 0.\]
Integrating by parts (in time) yields:
\begin{equation}
\label{eq:ADER_ibp}
\psi(\Delta t) u(\Delta t) - \psi(0) u(0)  
- \int^{\Delta t}_0 \partial_t \psi(t) u(t) {\rm d}t
+ \int^{\Delta t}_0 \psi(t) \partial_x f(u(t)) {\rm d}t = 0.
\end{equation}
Note that here we do not show the spatial dependency to simplify the notations. 
We now represent our solution using Lagrange polynomials {\it in time} $\ell_i(t)$ defined on $n+1$ Legendre quadrature points 
$\lbrace t_i \rbrace_{i=0}^n \in [0,\Delta t]$, which together with the quadrature weights $\lbrace w_i \rbrace_{i=0}^n$
can be used to perform integrals at the correct order in time. We are aiming at a solution with the same order of accuracy in time
than in space, so $n$ is taken here equal to the polynomial degree of the spatial discretisation.
We can write:
\[ u(t) = \sum_{i=0}^n u_i \ell_i(t),\]
and replace the integrals in Eq.~\eqref{eq:ADER_ibp} by the respective quadratures. We now replace the arbitrary test function $\psi(t)$
by the set of Lagrange polynomials $\{\ell_j(t)\}_{i=0}^n$ and obtain:

\begin{equation}\label{eq:System}
\ell_j(\Delta t)\left(\sum_{i=0}^{n} u_i \ell_i(\Delta t)\right) - \ell_j(0)u(0)
- \Delta t \sum_{i=0}^{n} w_i \ell^\prime_j(t_i) u_i 
+ \Delta t \sum_{i=0}^{n} w_i \ell_j(t_i) \partial_x f(u_i) 
=0 .
\end{equation}
To derive the previous equation, we used the interpolation property of the Lagrange polynomials with $u(t_i) = u_i$.
Note that $u(0)$ corresponds to the solution at the beginning of the time step.
The previous system can be rewritten in a matrix form, defining a mass matrix $M \in \mathbb{R}^{(n+1)\times(n+1)}$ 
and a right-hand side vector $r$ as:
\begin{equation}\label{eq:MassmatrixAder}
M_{ji} = \ell_j(\Delta t)\ell_i(\Delta t)- \Delta t w_i \ell^\prime_j(t_i)~~~{\rm and}~~~r_j = \ell_j(0)u(0) - \Delta t \sum_{i=0}^{n} w_i \ell_j(t_i) \partial_x f(u_i).
\end{equation}
The previous implicit non-linear equation with unknown $\lbrace u_i \rbrace_{i=0}^n$, is now written as:
\begin{equation}\label{fix:point}
 M_{ji} u_i = r_j(u_0,...,u_n),
\end{equation}
which can be solved with a fixed-point iteration method. We use a uniform initial guess with $\lbrace u^0_i=u(0)\rbrace_{i=0}^n$ and perform a standard Picard iterative scheme as follows
\begin{equation}
\label{eq:fixpoint_iteration}
u_i^{k+1}=M^{-1}_{ij} r_j (u_0^k,...,u_n^k),
\end{equation}
where index $k$ stands for the iteration count.
Finally, we use these final predicted states $\lbrace u^{n}_i \rbrace_{i=0}^n$ at our quadrature points and update the final solution as:
\begin{equation}
u(\Delta t) = u(0) - \Delta t \sum_{i=0}^n w_i \partial_x f(u^{n}_i).
\end{equation}

Because we always have this final update, we only need $n$ internal corrections to the solution (iterations) to obtain a solution that is accurate up to order $n+1$ in time \citep{dumbser_ader_2013}. The first order scheme with $n=0$ does not require any iteration, as it uses only the initial first guess to compute the final update, corresponding exactly to the first-order forward Euler scheme.

Note that in this flavour of the ADER scheme, we need to estimate the derivative of the flux for each time slice according to the SD method, including the Riemann solvers at element boundaries, making it different from the traditional ADER-DG framework presented in \cite{dumbser_ader_2013} and more similar with \cite{mhv2020}, which remains local until the final update. Precisely because we include the Riemann solver at the element boundaries, we maintain the continuity requirement on the normal component, needed for the appropriate evolution of the induction equation.
We use a Courant stability condition adapted to the SD scheme, as explained in \cite{vanharen2017} and compute the time step as:
\[
\Delta t = \frac{C}{n+1} \frac{\Delta x}{|v_{\rm max}|},
\]
where again $C=0.8$ and $n$ is the polynomial degree of our discretisation in space. We justify this choice by a careful time stability analysis in the following section.

\subsection{A modified SD scheme for the induction equation}

\begin{figure}
    \centering
    \includegraphics[width=1.0\textwidth]{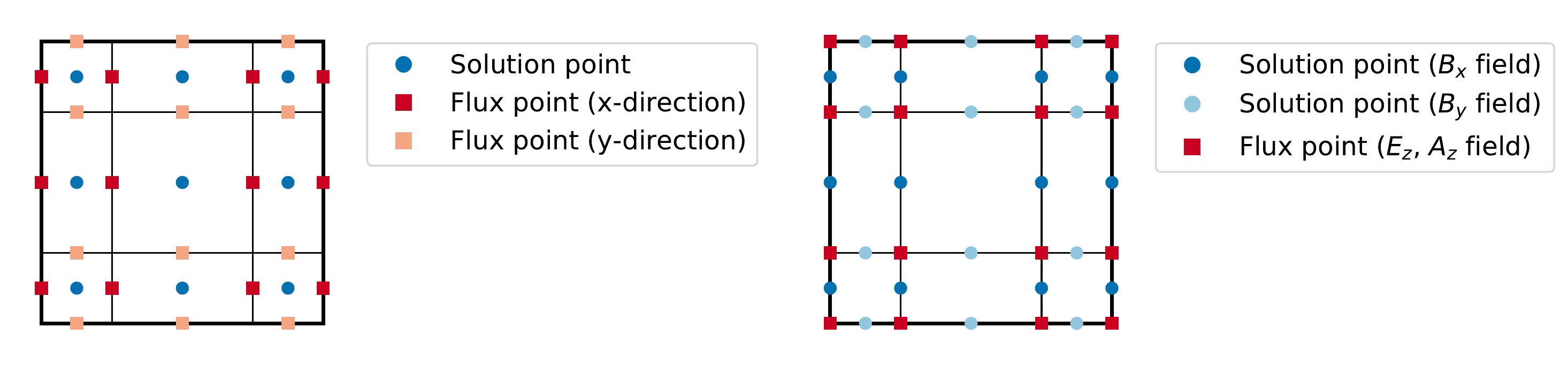}
    \caption{Left panel: position of the solution points in blue, and of the flux points in red, for a traditional SD method with $n=2$. 
    Right panel: position of the solution points for $B_x$ and $B_y$ in blue and of the flux points for the electric field $E_z$ and the vector potential $A_z$ in red, for our new SD method for the induction equation with $n=2$. \label{fig:sd_representation}}
\end{figure}

The traditional SD method is particularly well suited for the Euler sub-system with conservation laws based on the divergence operator.
In Fig.~\ref{fig:sd_representation}, we show on the left panel the traditional discretisation of one element using SD, with the control volume boundaries shown 
as black solid lines, the solution points in blue inside the control volumes, and the flux points in red on the faces of each control volume.
The strict conservation property of SD can be explained using for example the density. Defining the corner points of the control volumes as $( x_i, y_j)$,
we can compute the total mass within a rectangle defined by the four points $( 0, 0)$, $( x_i, 0)$, $( 0, y_j)$ and $( x_i, y_j)$ as
$M(x_i,y_j)$. Note that the corner points are defined as the intersection of the lines where the flux points are defined. 

We now represent this cumulative mass everywhere inside the element using Lagrange polynomials defined on the flux points as
\begin{equation}
    M(x,y) = \sum_{i=0}^{n+1} \sum_{j=0}^{n+1} M(x_i,y_j)\ell_i(x) \ell_j(y),
\end{equation}
where we dropped the superscript $f$ for simplicity. The density field is obtained by taking the derivative of the cumulative mass as:
\begin{equation}
    \rho(x,y) = \sum_{i=0}^{n+1} \sum_{j=0}^{n+1} M(x_i,y_j)\ell^{\prime}_i(x) \ell^{\prime}_j(y),
\end{equation}
where the prime stands for the spatial derivative. This exact procedure can be used to initialise the value of $\rho$ at the solution points, 
as well as to prove that the SD method as described in the previous section 
is strictly conservative \cite{Liu2004,Liu2006}. 

The induction equation, however, is a conservation law for the magnetic flux through a surface, as it does not feature a divergence operator
but a curl operator. We therefore propose a small modification of the classical SD scheme, similar to that of \citep{chandrashekar_2020,praveen_2019}, with different collocation points for the magnetic field components $B_x$ and $B_y$.

In the two-dimensional case, we start with a vector potential $\vec{A} = (0,0,A_z)$. We approximate the z-component using the red square collocation points, as denoted in the right panel of Fig. \ref{fig:sd_representation}. The approximation takes the form

\[
A_z(x,y) = \sum_{i=0}^{n+1}\sum_{j=0}^{n+1} A_z(x_i,y_j)\ell_i(x)\ell_j(y)
\]

Then, the magnetic field, in 2-dimensions is obtained by:

\[
B_x(x,y) = \sum_{i=0}^{n+1}\sum_{j=0}^{n+1} A_z(x_i,y_j)\ell_i(x)\ell^{'}_j(y) \quad 
B_y(x,y) = -\sum_{i=0}^{n+1}\sum_{j=0}^{n+1} A_z(x_i,y_j)\ell^{'}_i(x)\ell_j(y)
\]

Because of how the magnetic field $\vec{B} = (B_x, B_y)$ is initialised, this is by definition divergence-free.

Then, we define the magnetic flux $\phi_x$ through the surface defined by two corner points $(x_i,0)$ and $(x_i,y_j)$ as:
\begin{equation}
\label{eq:mfx}
\phi_x(x_i,y_i) = \int_0^{y_i} B_x(x_i,y) {\rm d} y.
\end{equation} 
Similarly, we define the magnetic flux $\phi_y$ through the surface defined by two corner points $(0,y_j)$ and $(x_i,y_j)$ as:
\begin{equation}
\label{eq:mfy}
\phi_y(x_i,y_i) = \int_0^{x_i} B_y(x,y_i) {\rm d}x.
\end{equation}
We see that $\phi_x$ and $\phi_y$ are both defined over the set of corner points $(x_i,y_j)$.
We can now represent the numerical approximation $\phi_x$ (resp. $\phi_y$) using Lagrange polynomials defined on the flux points as:
\[\phi_x(x,y) = \sum_{i=0}^{n+1} \sum_{j=0}^{n+1} \phi_x(x_i,y_j) \ell_i(x) \ell_j(y). \] 
Then, we deduce the numerical approximation of $B_{x,h}$ as:
\begin{equation}
\label{eq:sd_bx}
B_{x}(x,y) = \partial_y \phi_x = \sum_{i=0}^{n+1} \sum_{j=0}^{n+1} \phi_x(x_i,y_j) \ell_i(x) \ell^{\prime}_j(y),
\end{equation}
and the numerical approximation of $B_{y,h}$ as:
\begin{equation}
\label{eq:sd_by}
B_{y}(x,y) = \partial_x \phi_y = \sum_{i=0}^{n+1} \sum_{j=0}^{n+1} \phi_y(x_i,y_j) \ell^{\prime}_i(x) \ell_j(y).  
\end{equation}
The key difference between this configuration and the traditional SD method is that for $B_x$, the $x$ direction has an extra degree of freedom 
and a higher polynomial degree (similarly for $B_y$ in the y direction). 
This also means that the corresponding solution points for $B_x$ and $B_y$ are staggered with respect to the traditional SD method. 
In the right panel of Fig~\ref{fig:sd_representation}, we show the position of these new solution points for $B_x$,  $B_y$ (in blue) 
and new flux points (in red) where the electric field will be defined, as explained below. 
Note that if the initial magnetic field is divergence free, applying the divergence theorem to the rectangle defined by the same four corner points as before leads to
the constraint:
\begin{equation}
\label{eq:discrete_circulation}
\phi_x(x_i,y_j) + \phi_y(x_i,y_j) - \phi_x(0,y_j) - \phi_y(x_i,0) = 0, \quad \forall i, j. 
\end{equation}

\begin{proposition}
\label{proposition:eq_25}
Equation \eqref{eq:discrete_circulation} holds if we can integrate $\vec{B}$ exactly or by starting from a vector potential $\vec{A}$.
\end{proposition}
\begin{proof}
Using the numerical approximation of $A_z$ (and sub-consequently of $\vec{B}$), we can write the fields $\phi_x(x,y)$ and $\phi_y(x,y)$ for any control volume $K = [0,x_m]\times [0,y_m]$

\begin{align*}
\phi_x(x_m,y_m) &= \int_0^{y_m} \sum_{i=0}^{n+1}\sum_{j=0}^{n+1} A_z(x_i,y_j)\ell_i(x_m)\ell^{'}_j(y) dy \\ &= \sum_{i=0}^{n+1}\sum_{j=0}^{n+1} A_z(x_i,y_j)\ell_i(x_m) \int_0^{y_m} \ell^{'}_j(y) dy \\
&= \sum_{i=0}^{n+1}\sum_{j=0}^{n+1} A_z(x_i,y_j)\ell_i(x_m) \left[\ell_j(y_m) - \ell_j(0) \right] 
\end{align*}

\begin{align*}
\phi_y(x_m,y_m) &= -\int_0^{x_m}  \sum_{i=0}^{n+1}\sum_{j=0}^{n+1} A_z(x_i,y_j)\ell_j(y_m)\ell^{'}_i(x) dx \\ &= - \sum_{i=0}^{n+1}\sum_{j=0}^{n+1} A_z(x_i,y_j)\ell_j(y_m) \int_0^{x_m}\ell^{'}_i(x) dx \\ &= - \sum_{i=0}^{n+1}\sum_{j=0}^{n+1} A_z(x_i,y_j)\ell_j(y_m) \left[\ell_i(x_m) - \ell_i(0) \right]
\end{align*}

The integration is exact given the right quadrature rule. 

Then, we can observe that \eqref{eq:discrete_circulation} holds.

\end{proof}

\begin{proposition}
\label{proposition:pointwise_div_free}
The proposed numerical representation of $\vec{B}$ is pointwise (or locally) strictly divergence free.
\end{proposition}
\begin{proof}
We now evaluate the divergence of the numerical approximation $\vec{B}(x,y) = [B_x,B_y]$:
\begingroup
\allowdisplaybreaks
\begin{align*}
\partial_x B_x + \partial_y B_y &= 
\partial_x \left(\sum_{i=0}^{n+1} \sum_{j=0}^{n+1} \phi_x(x_i,y_j) \ell_i(x) \ell^\prime_j(y)\right) + 
\partial_y \left(\sum_{i=0}^{n+1} \sum_{j=0}^{n+1} \phi_y(x_i,y_j) \ell^\prime_i(x) \ell_j(y)\right)\\
&= \sum_{i=0}^{n+1} \sum_{j=0}^{n+1} \left( \phi_x(x_i,y_j) + \phi_y(x_i,y_j) \right)\ell_i'(x) \ell_j'(y)\\
&= \sum_{i=0}^{n+1} \sum_{j=0}^{n+1} \left( \phi_x(0,y_j) + \phi_y(x_i,0) \right)\ell_i'(x) \ell_j'(y),
\end{align*}
\endgroup
where we used the property that the total magnetic flux through the rectangle vanishes (see Eq.~\eqref{eq:discrete_circulation}).
We can now separate and factor out the $i$ and $j$ sums as:
\begingroup
\allowdisplaybreaks
\begin{align*}
\partial_x B_x + \partial_y B_y
&= \sum_{i=0}^{n+1} \sum_{j=0}^{n+1}\phi_x(0,y_j)\ell_i'(x) \ell_j'(y) + \sum_{i=0}^{n+1} \sum_{j=0}^{n+1} \phi_y(x_i,0) \ell_i'(x) \ell_j'(y)\\
&= \left( \sum_{j=0}^{n+1}\phi_x(0,y_j) \ell_j'(y)\right) \left(\sum_{i=0}^{n+1}\ell_i'(x)\right) + \left(\sum_{i=0}^{n+1} \phi_y(x_i,0) \ell_i'(x)\right)\left(\sum_{j=0}^{n+1} \ell_j'(y)\right)\\
&= 0,
\end{align*}
\endgroup
where we used the property of the Lagrange polynomials that $\sum_{i=0}^{n+1} \ell_i (x)=1$ so that the corresponding derivative vanishes uniformly.
\end{proof}

\begin{proposition}
\label{proposition:globally_div_free}
The proposed numerical representation of $\vec{B}$ is globally divergence free.
\end{proposition}
\begin{proof}
If the initial magnetic field is divergence free, $B_x$ is continuous across the left and right boundaries of each element. Similarly, $B_y$ is continuous across the bottom and top boundaries of the element. It follows that $\phi_x$ (resp. $\phi_y$) is initially identical on the left (resp. bottom) edge of the right (resp. top) element and on the right (resp. top) edge of the left (resp. bottom) element. Because the adopted Lagrange polynomial basis is an interpolatory basis, and because the solution points of $B_x$ and $B_y$ are collocated on the element boundaries, the continuity of the magnetic field in the component normal to the element face is enforced by construction and at all orders. Note that the case $n=0$ corresponds exactly to the Constrained Transport method, as implemented in popular FV codes. The proposed discretisation is a generalisation of CT to arbitrary high order.
\end{proof}

We now describe the SD update for the induction equation. We define the electric field $\vec{E}= - \vec{v}\times\vec{B}$ and write the induction equation as
\[\partial_t \vec{B} = - \nabla\times \vec{E} . \]
Once we know the prescribed velocity field and the polynomial representation of the magnetic field throughout the element,
as in Eq.~\eqref{eq:sd_bx} and Eq.~\eqref{eq:sd_by}, we can compute the electric field at the control volume corner points $(x_i,y_j)$. These are the equivalent of the flux points in the traditional SD method.
Since the electric field is continuous across element boundaries, we need to use a 1D Riemann solver for flux points inside the element edges, and a 2D Riemann solver
at corner points between elements. This step is crucial as it maintains the global divergence-free property. We see for example that the electric field on an element face
will be identical to the electric field on the same face of a neighbouring element. 2D Riemann solvers at element corners are also important to maintain this global property, and in the case of the induction equation, we just need to determine the 2D upwind direction using both $v_x$ and $v_y$. After we have enforced a single value for the electric field on element
edges, we can interpolate the electric field inside the element, using flux points and the corresponding Lagrange polynomials, as before:
\begin{equation}
E_z(x,y) = \sum_{i=0}^{n+1} \sum_{j=0}^{n+1} E_z(x_i,y_j) \ell_i(x) \ell_j(y),
\end{equation}
and update the magnetic field directly using the pointwise update:
\begin{equation}
\label{eq:induction_update}
\partial_t B_{x} = - \partial_y E_z,~~~{\rm and}~~~ \partial_t B_{y} = \partial_x E_z.
\end{equation}
We have directly:
\begin{equation}
\partial_t \left( \partial_x B_{x} + \partial_y B_{y} \right) = 0,
\end{equation}
so that the divergence of the field, if zero initially, will remain zero at all time.
The continuity of $E_z$ at element boundaries also implies that the continuity of the normal components of the magnetic field will be preserved after the update.

Note that at the beginning of the time step, we only need to know the values of $B_x$ and $B_y$ at their corresponding solution points to obtain the same polynomial interpolation  
as the one we derived using the magnetic fluxes. This follows from the uniqueness of the representation of the solution by polynomials of degree $n$. Similarly, the time update we just described
can be performed only for the magnetic solution points (see Fig.~\ref{fig:sd_representation}) to fully specify our zero divergence field for the next time step.

\begin{algorithm}[H]
 \KwData{ $A_z$ at $t=0$}
 \KwResult{$\vec{B}$ at $t=T$ }
 compute the numerical representation $A_z$\;
 build $\phi_x$ and $\phi_y$ by integrating $B_x$ and $B_y$ (which are given by differentiating $A_z$)\;
 get $B_x$ and $B_y$ by differentiating $\phi_x$ and $\phi_y$\;
 \While{t < T}{
    perform ADER-SD update on nodal values of $B_x$ and $B_y$ through \eqref{eq:induction_update}\;
    %evolve nodal values of $B_x^h$ and $B_y^h$ placed as in figure \ref{fig:sd_representation} as \eqref{eq:induction_update} \;
 }
 \caption{SD-ADER algorithm compatible with the induction equation.}
\end{algorithm}

\begin{algorithm}[H]
 \KwData{$B_x$ and $B_y$  at $t=t^n$}
 \KwResult{$B_x$ and $B_y$  at $t=t^{n+1}$}
 \While{iteration < total iterations}{
     compute $E_z$ field on flux points (refer to Fig.  \ref{fig:sd_representation}) for all time-substeps\;
     compute unique value of $E_z$ using a 1-dimensional Riemann solver at cell faces and a 2-dimensional Riemann solver at cell corner points\;
     build $E_z$ flux in space-time\;
     perform ADER sub-timestep update on degrees of freedom of $B_x$ and $B_y$ \;
 }
 \caption{ADER-SD update}
\end{algorithm}

\begin{remark}
The idea of approximating the magnetic field $\vec{B}$ through tensor product polynomials while keeping $\vec{B}\cdot\vec{n}$ continuous across cell faces is a well known idea, for example, through the use of Raviart-Thomas (RT) elements \cite{Brezzi_1991} in the finite element context. In fact, this approach has been used to treat the induction equation \cite{Balsara_weno_2009,praveen_2019}. The main difference between our method and RT is that we do not explicitly have to build RT approximation basis. In particular, the continuity of  $\vec{B}\cdot\vec{n}$ across cells is guaranteed by exact interpolation of nodal values collocated appropriately, as well as a Constrained Transport-like update using a unique electric field $E_z$.
\end{remark}

\begin{proposition}
The previous scheme is equivalent to a simple evolution equation for the vector potential with a continuous SD scheme, for which both the solution points and the flux points are defined on the corner points of the control volumes.
\end{proposition}
\begin{proof}
The magnetic fluxes introduced in \eqref{eq:mfx} and \eqref{eq:mfy} are analogous to a magnetic potential in two space dimensions. 
Indeed, in this case, one can compute directly the magnetic vector potential $\vec{A} = (0,0,A_z)$ at the corner points, using 
\begin{equation}
A_z(x_i,y_i) = \phi_x(x_i,y_i) - \phi_y(x_i,0) = \phi_x(0,y_i) - \phi_y(x_i,y_i).
\end{equation}
We then interpolate the vector potential within the elements using Lagrange polynomials as:
\begin{equation}
A_z(x,y) = \sum_{i=0}^{n+1} \sum_{j=0}^{n+1} A_z(x_i,y_j) \ell_i(x) \ell_j(y)
\end{equation}
and compute the magnetic field components as:
\begin{equation}
B_{x}(x,y) = \partial_y A_z~~~{\rm and}~~~ B_{y}(x,y) = - \partial_x A_z.
\end{equation}
This definition is equivalent to the previous one. For the SD update, we compute the electric field at each corner point, using again a Riemann solver for multi-valued flux points.
The vector potential is then updated directly at the corner point using
\begin{equation}
\label{eq:HJ}
\partial_t A_z = - E_z = v_x B_y - v_y B_x = -v_x \partial_x A_z - v_y \partial_y A_z.
\end{equation}
We see that this last equation yields an evolution equation for $A_z$, where all terms are evaluated at the corner points. This corresponds to a variant of the SD method, for which the solution points are not placed at the centre of the control volumes, but migrated for example to their upper-right corner, and for which the flux points are not placed at the centre of the faces, but migrated to the same upper-right corner (thus overlapping). Note however an important difference with the traditional SD method: the vector potential $A_z$ is a continuous function of both $x$ and $y$ so that we have $n+2$ solution points instead of $n+1$. In other words, each element face shares the same values for $A_z$ and $E_z$ with the corresponding neighbouring element face. We have therefore a strict equivalence between the induction equation solved using our SD scheme for the magnetic field and the evolution equation solved using this particular SD variant for the vector potential. 
\end{proof}

\subsection{Stability for the linear induction equation}
\label{sec:stability}

\begin{figure}
\centering
\includegraphics[width=.48\textwidth]{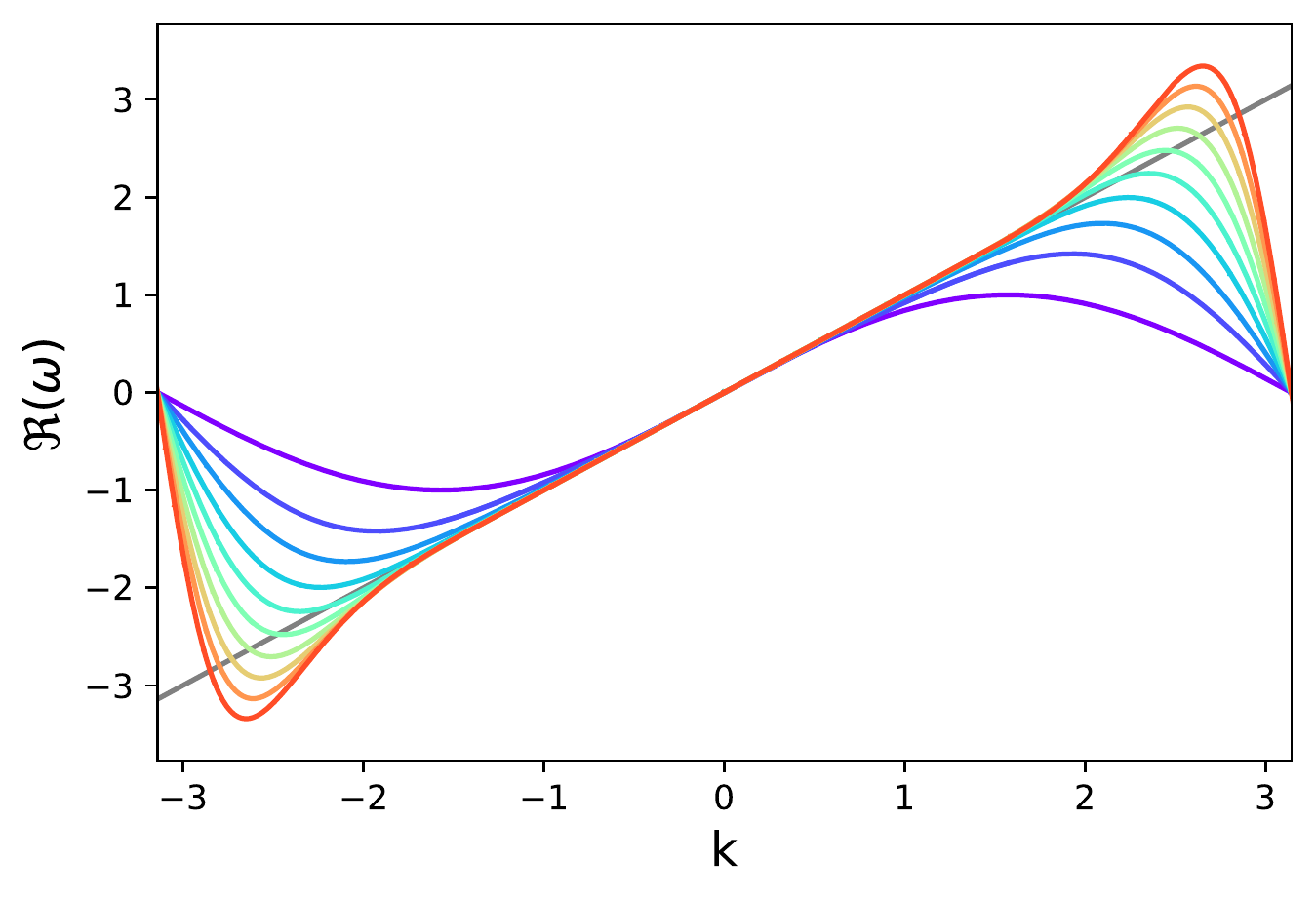}
\includegraphics[width=.48\textwidth]{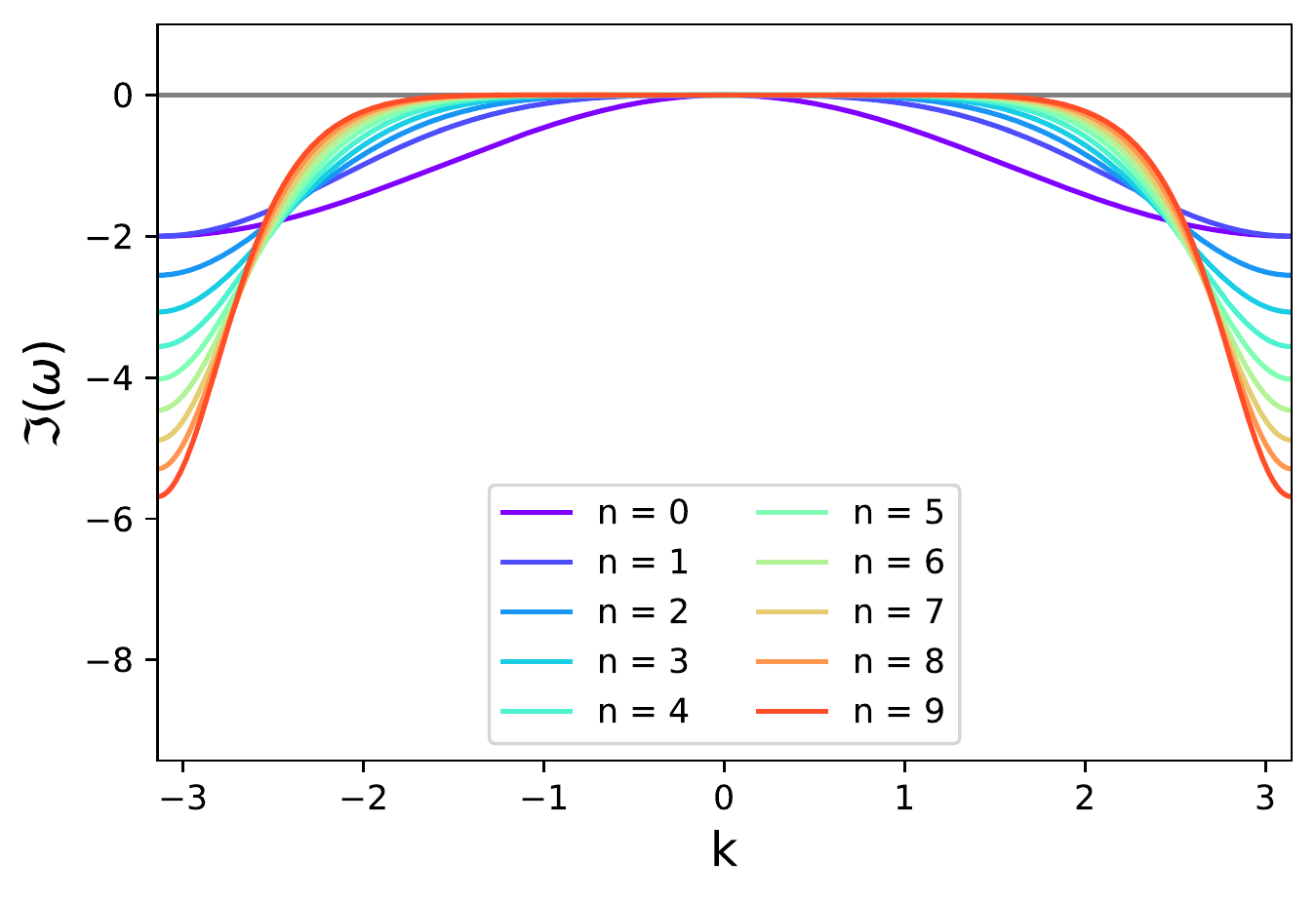}
\caption{Real (left panel) and imaginary (right panel) parts of $\omega$ for different polynomial degrees $n$. In each panel, the grey line represents the exact dispersion (left) and diffusion (right) relation for the waves. The wave number is expressed here in units of $(n+1)/\Delta x$, while the frequency $\omega$ is shown in units of $(n+1) v /\Delta x$. \label{fig:stability}}
\end{figure}

We will now demonstrate that the proposed SD scheme for the induction equation is stable. We will also analyze the dispersion properties of the scheme at various orders.
To achieve this goal, we will first exploit the strict equivalence between the magnetic field update and the vector potential update, as explained previously.
We will prove the stability of the scheme in one space dimension, and use the tensor product rule to extend these results to higher dimensions. 

We write our vector potential equation in 1D, assuming here, without loss of generality, a positive velocity field $v>0$: 
\begin{equation}
\label{eq:vector_potential_pde}
    \partial_t A + v \partial_x A = 0.
\end{equation}
Space is discretised using N equal-size element with width $\Delta x = L/N$ and labelled by a superscript $p=1,\dots,N$. We have $n+2$ solution points for the vector potential,
which coincide with the flux points of the SD method, here labelled as $x_i^p$ with $i=0,\dots,n+1$. As explained earlier, we use for the inner $n$ flux points 
the Gauss-Legendre quadrature points, while the leftmost one, $x_0^p$, is aligned with the left element boundary, and the rightmost one, $x_{n+1}^p$,
is aligned with the right element boundary. We see that we have redundant information in our solution vector $A_i^p$ as $x_{0}^p = x_{n+1}^{p-1}$,
so that $A_0^p = A_{n+1}^p$ and $A_{n+1}^p=A_0^{p+1}$. This redundancy is a fundamental difference with the traditional SD method and ensures that the vector 
potential is continuous across element boundaries. In our present derivation, we need to avoid this duplication and define the solution vector $A_i^p$
in each element only for $i=1,\dots,n+1$, dropping the leftmost point and assigning it to the left element. This choice is arbitrary, but it corresponds here 
to the upwind solution of the Riemann solver at the left boundary, as we have $v>0$. Our vector potential solution points now resemble the classical SD method solution 
points shifted to the right of their closest flux points. The vector potential is interpolated within element $p$ using the $n+2$ flux points as:
\begin{equation}
    A^p(x) = A_{n+1}^{p-1} \ell_0(x) + \sum_{j=1}^{n+1} A_j^p \ell_j(x).
\end{equation}
We can write the corresponding SD update as
\begin{equation}
    \partial_t A_i^p = -v \left( A_{n+1}^{p-1} \ell^\prime_0(x_i) + \sum_{j=1}^{n+1} A_j^p \ell^\prime_j(x_i) \right).
\end{equation}

For the stability analysis, we follow the methodology presented in \cite{hu1999, abeele2008} and study the response of the scheme to a planar wave solution of the form:
\begin{equation}
    A(x) = \tilde A \exp(i ( k x -\omega t )),
\end{equation}
using periodic boundary conditions. The stability of a planar wave solution will depend on the imaginary part of $\omega$. Indeed, the amplitude of the wave will not increase if $\Im(\omega)$ remains smaller than 0.

The flux points coordinates are split between the element leftmost coordinates and the relative flux point coordinates as
$x_i^p = (p-1) \Delta x + x_i$, so that we have:
\begin{equation}
A_i^p = \tilde A \exp(-i \omega t)\exp(i k (p-1) \Delta x) \exp(i k x_i).
\end{equation}
The update now writes
\begin{equation}
-i \omega \tilde A \exp(i k x_i) = -v \left( \tilde A \exp(i k x_{n+1}) \ell^\prime_0(x_i) \exp(-i k \Delta x) + \sum_{j=1}^{n+1} \tilde A \exp(i k x_{j}) \ell^\prime_j(x_i)  \right) .
\end{equation}
We define the solution vector for the planar wave as $u_i = \tilde A \exp(i k x_i)$. The previous equation can be written in matrix form as follows:
\begin{equation}
\left( -i \frac{\omega \Delta x}{v} \mathbb{I} + \mathbb{M} \right) \vec{u} = 0, 
\end{equation}
where, for sake of simplicity, we have considered a normalized coordinate system inside each element so that $x_0=0$ and $x_{n+1}=1$. This explains why
the factor $\Delta x$ has been factored out. The matrix $\mathbb{M}$ is defined as:
\begin{equation}
\mathbb{M} =
\begin{bmatrix}
\ell^\prime_1(x_1) & \ldots & \ell^\prime_n(x_1) & \ell^\prime_{n+1}(x_1) + \ell^\prime_0(x_1) \exp(-i k \Delta x) \\
\ldots & \ldots & \ldots & \ldots\\
\ell^\prime_1(x_{n+1}) & \ldots & \ell^\prime_n(x_{n+1}) & \ell^\prime_{n+1}(x_{n+1}) + \ell^\prime_0(x_{n+1}) \exp(-i k \Delta x)
\end{bmatrix}.
\end{equation}
The dispersion relation of the waves is obtained by requiring 
\begin{equation}
\det \left( -i \frac{\omega \Delta x}{v} \mathbb{I} + \mathbb{M} \right) = 0 ,
\end{equation}
%which is equivalent to finding the $n+1$ complex eigenvalues of matrix $\mathbb{M}$.
which amounts to finding the $n+1$ complex eigenvalues of matrix $-i\mathbb{M}$. We can then represent the dispersion relation of the scheme with
$\Re{(\omega)}$ and the diffusion relation with $\Im{(\omega)}$. More importantly, the wave amplitude will be damped if $\Im{(\omega)}<0$, corresponding to a stable numerical scheme, and will be amplified exponentially if $\Im{(\omega)}>0$, corresponding to an unstable numerical scheme. The maximum wave number is set by the grid size $\Delta x$
and the polynomial degree $n$ so that:
\begin{equation}
    k_{\rm max} = \left( n+1 \right) \frac{\pi}{\Delta x} = \left( n+1 \right) N \frac{\pi}{L}.
\end{equation}
We see that the maximum wave number depends on the product $(n+1)\times N$ which corresponds to the number of degrees of freedom of the SD method.
The previous dispersion relation generates $n+1$ eigenvalues in the k-interval $\left[ -\pi/\Delta x, \pi/\Delta x\right]$, 
owing to the periodicity of the function $\exp(-i k \Delta x)$. In order to derive the dispersion relation in the entire range of wave number 
$\left[ -(n+1)\pi/\Delta x, (n+1)\pi/\Delta x\right]$, 
the eigenvalues have to be shifted by an integer multiple of $2\pi/\Delta x$ to generate a single branch in the dispersion relation.

We show in Fig.~\ref{fig:stability} the real and imaginary part of $\omega$ for a set of SD schemes that have exactly the same number of degrees of freedom $(n+1) \times N$,
with $n$ ranging from 0 to 9.
We note that,
although our scheme is different from the classical SD scheme, the dispersion relations at these various orders are identical to the corresponding dispersion relation found by \cite{abeele2008} for the classical SD method. This strict equivalence is true only for a constant velocity field. We see also in Fig.~\ref{fig:stability} that, although all these
schemes have exactly the same number of degrees of freedom, the higher the polynomial degree, the closer the scheme gets to the true solution, 
namely $\Im{(\omega)}=0$ and $\Re{(\omega)}=v k$, shown as a grey line in Fig.~\ref{fig:stability}.
We conclude from this analysis that the SD spatial operator is stable, because $\Im{(\omega)}<0$ everywhere. To explicitly connect these results to \cite{abeele2008}, one can see that the Fourier footprint $\Omega$ can be obtained from the relation $\Omega = -i\omega$. With this nomenclature, our scheme has $\mathcal{R}(\Omega) < 0$.
%\end{comment}

The SD semi-discretisation of the PDE \eqref{eq:vector_potential_pde} leads to a system of first order ordinary differential equations in time:
\begin{equation}
    \begin{cases}
          U'(t) &= F(U) \\
          U(0) &= U_0.
    \end{cases}
\end{equation}
We note $DOF$ the total number of degrees of freedom of the semi-discrete spatial SD operator defined by $U(t):\mathbb{R}\to\mathbb{R}^{DOF}$ and $F:\mathbb{R}^{DOF}\to\mathbb{R}^{DOF}$, respectively the vectors of unknowns and the discrete operator in space for all the degrees of freedom. We now show that using the ADER time stepping strategy, we  obtain a stable, fully discrete in space and time, high-order numerical scheme.

\begin{figure}
\centering
\includegraphics[width=.55\textwidth]{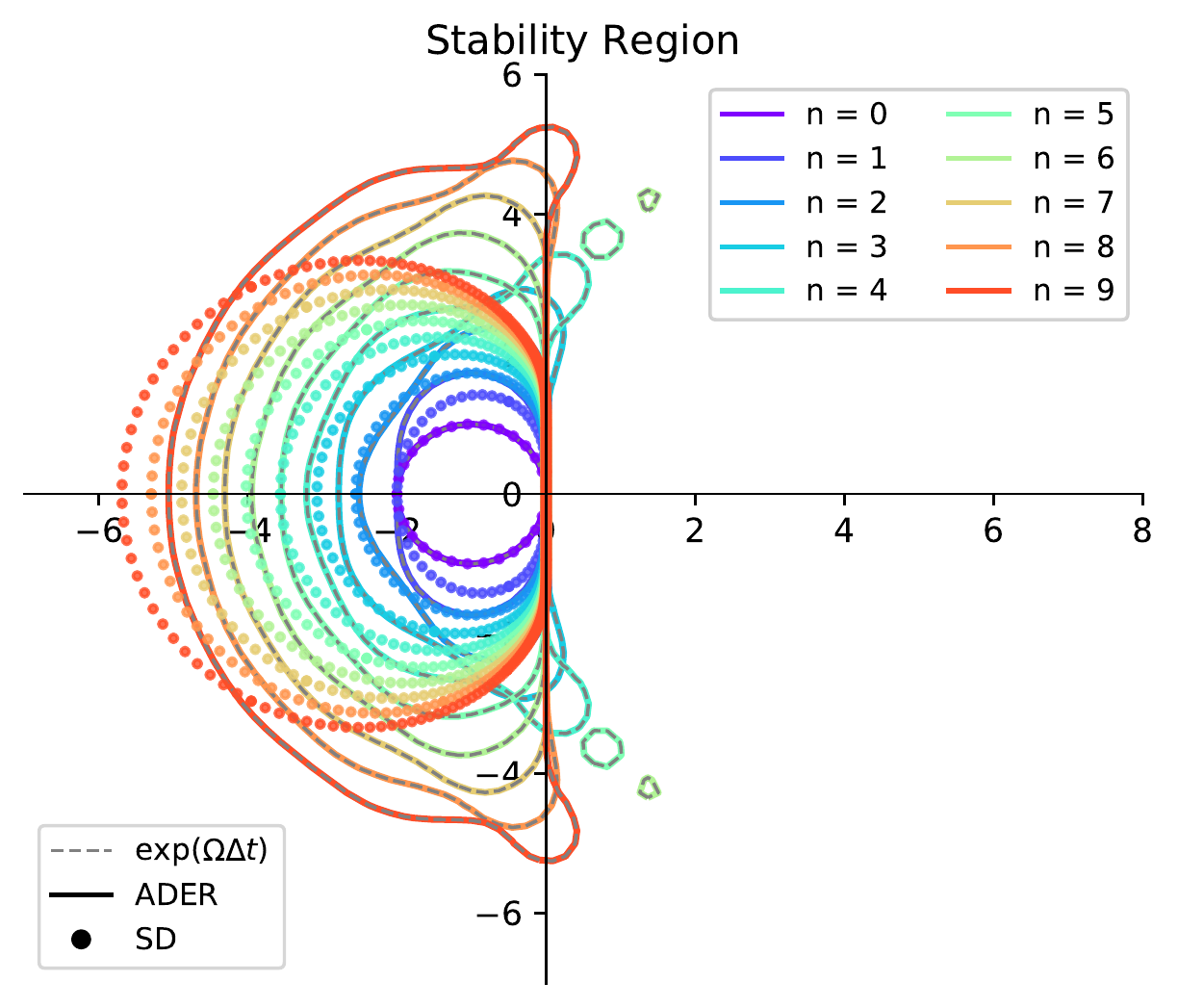}
\caption{Stability limits for ADER methods in the complex $\Omega \Delta t$-plane (continuous lines), from 0 to 9 corrections (ADER0 to ADER9), together with the stability domains of the SD space discretisation (symbols) for CFL = 1.0. Note that the stability region of the ADER scheme is identical to the exact amplification factor of $\exp(\Omega \Delta t)$. See text for details.
\label{fig:ader_stability}}
\end{figure}

We can investigate the full system of semi-discretised equations by isolating a single mode. Taking an eigenvalue $\Omega$ of the spatial discretisation operator,  we consider the canonical ODE:

\[ \frac{d}{dt} u = \Omega u.\]
We can write a general time integration method as
\[u^{n+1} = P(\Omega \Delta t) \cdot u^n,\]
where the operator $P$, called the numerical amplification factor, depends on the single parameter $\Omega \Delta t$. If we designate the eigenvalues of $P$ as $z_P$, the necessary stability condition is that all the eigenvalues $z_p$ should be of modulus lower than, or equal to, one \cite{Hirsch1988NumericalCO}. Similarly to \cite{mhv2020}, we perform a numerical stability study of the ADER scheme presented in this paper. In Fig.~\ref{fig:ader_stability}, we show the stability domains of the SD scheme, together with ADER in the $\Omega\Delta t-$plane. We can note that from $n=2$ onwards, the CFL should be reduced to a value slightly smaller than unity. We note that the stability region that we obtain is the same as the exact amplification factor $\exp(\Omega \Delta t)$ up to order $n$. This is no surprise as all methods with $n$ stages (in our case, corrections) and order $n$ have the same domain of stability \cite{Hirsch1988NumericalCO}. Then, by choosing an appropriate CFL condition, we are able to guarantee that $z_P(\Omega \Delta t)$ remain inside the ADER stability region.

In the future, we would like to study the stability of our method in more detail, similarly to \cite{glaubitz2018application}, and, given the similarities between our work and the one presented in \cite{balsara_kappeli_2018}, a more detailed numerical study of the stability of this scheme is of high interest as well.
\section{Numerical results}
\label{sec:numerics}

\begin{figure}
\centering
\includegraphics[width=.48\textwidth]{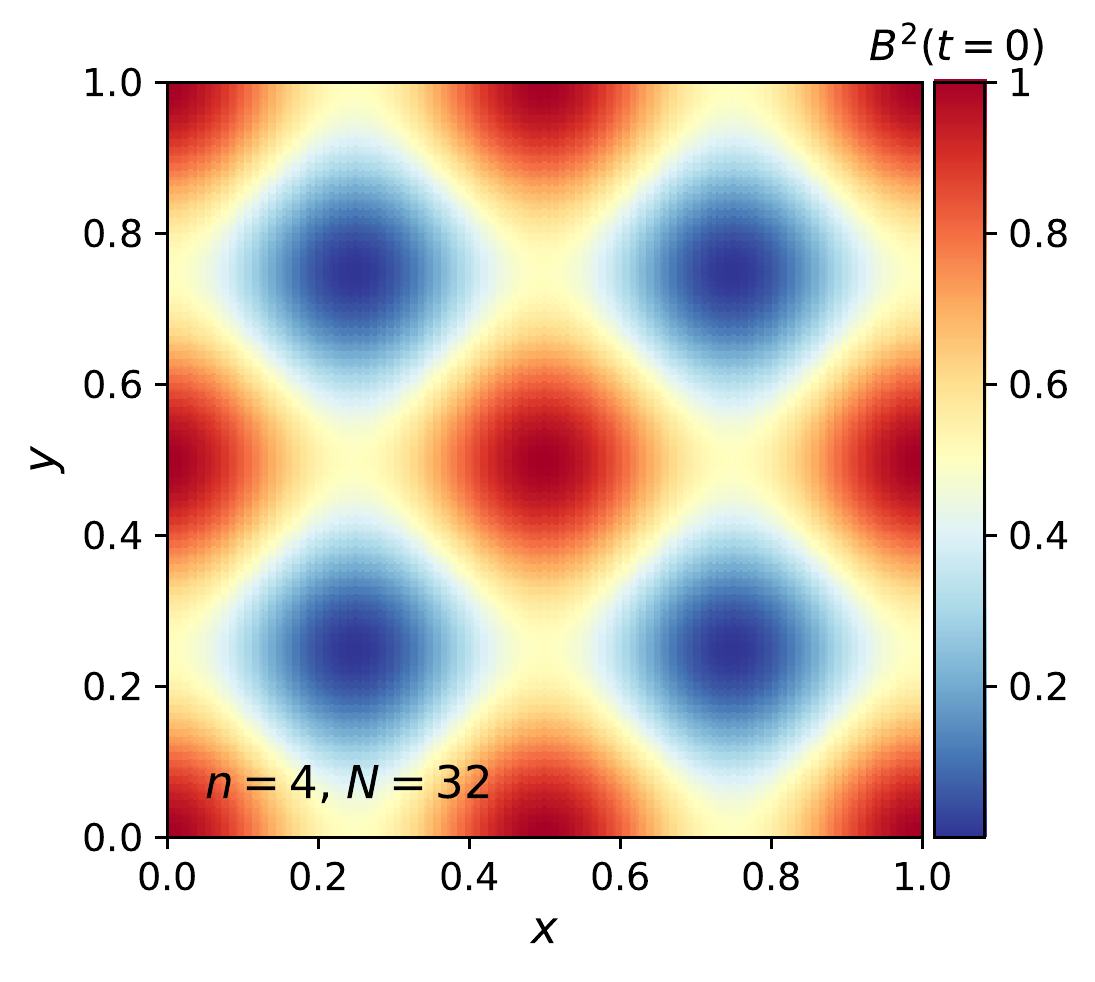}
\includegraphics[width=.5\textwidth]{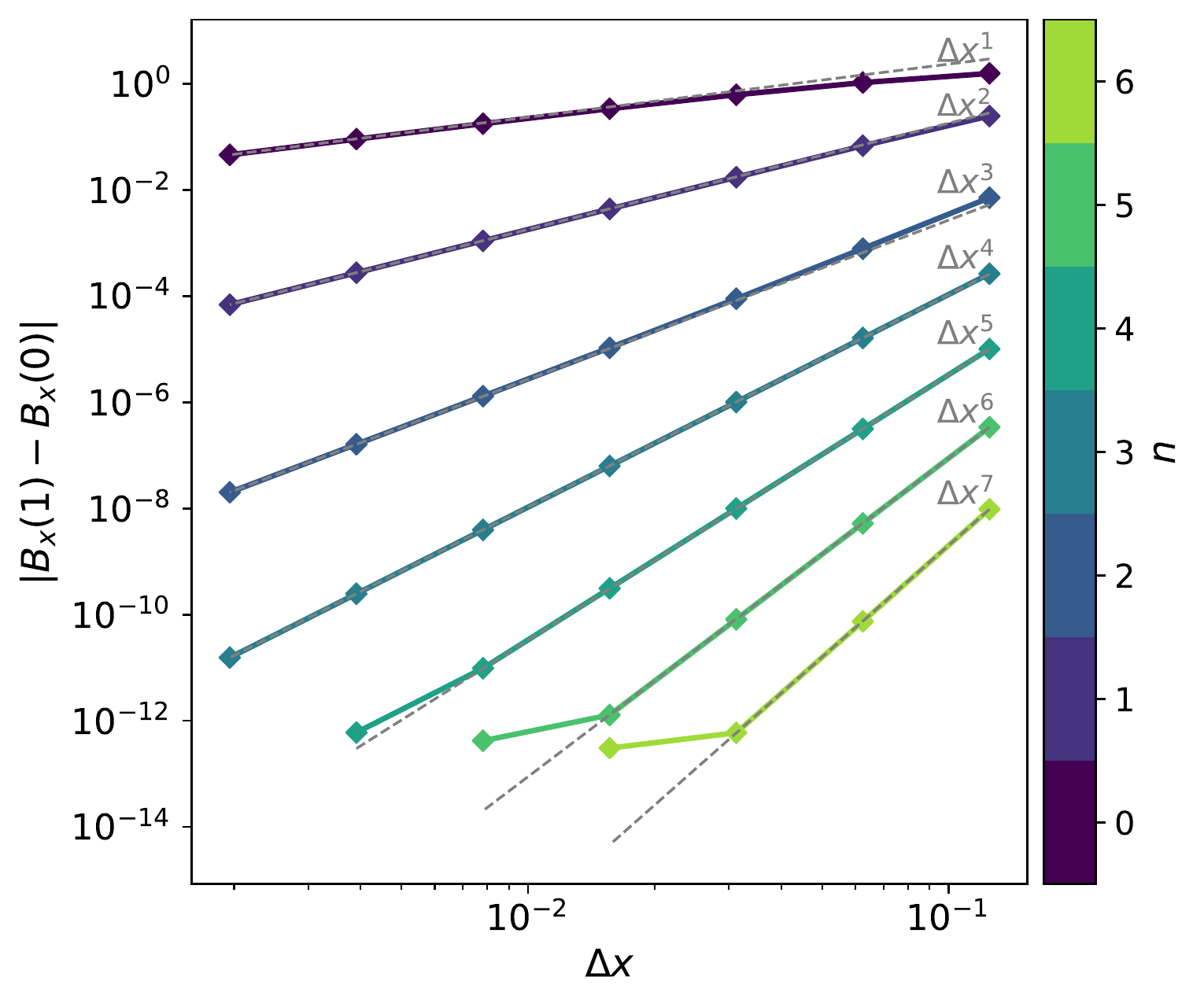}
\caption{Left panel: map of the magnetic energy density for the case $n=4$ and $N=32$. Right panel: $L^1$ convergence of the SD method for the smooth magnetic potential from Eq.~\eqref{eq:smoothloop} at different orders and spatial resolutions.
\label{fig:conv-rate-smooth-pot}}
\end{figure}

In this section, we test our new SD-ADER scheme for the induction equation using first a smooth initial condition, ensuring that our method is truly high order, and then using a more difficult tests, namely the advection of a discontinuous field loop under a constant velocity and a rotating velocity field. Finally, we compare our new SD-ADER scheme's performance to the various variants of the RKDG scheme on the advection of a discontinuous field loop problem.

\subsection{Continuous magnetic field loop}

In order to check that we are indeed solving the induction equation at the expected order of accuracy, we consider the advection of a smooth and periodic magnetic field given by the following initial conditions:
\begin{equation}
\label{eq:smoothloop}
 \vec{B} = \left( \cos(2\pi y), -\cos(2\pi x), 0\right),
\end{equation}
with a constant velocity field $\vec{v} = (1,1)$. 
We estimate the convergence rate of the proposed SD method by computing the $L^1$ error of each magnetic field component, averaged over the control volumes within each element. The $L^1$ error is defined as the $L^1$ norm of the difference in the mean value for each control volume between the numerical solution $u(t)$ at $t=1$ and the initial numerical approximation $u(0)$:
\begin{equation}
\label{eq:L1}
 L^1 = ||u(t)-u(0)||_1 = \sum_{K\in \mathcal{K}} \int_K |u(t)-u(0)|{\rm d} x{\rm d} y.
\end{equation}
In Fig.~\ref{fig:conv-rate-smooth-pot}, we present the $L^1$ convergence rates for $B_x$ only. We omit the results for $B_y$ as these are identical to the ones of $B_x$ due to the symmetry of the initial conditions. We can observe that the convergence rate of the method scales as $\Delta x^{n+1}$ (where $n$ is the polynomial degree of the interpolation Lagrange polynomial), as expected of a high-order method, and as observed in other high-order method implementations \cite{Schaal2015,Guillet2019,Derigs2018}.

As introduced in the previous section, the product $(n+1)
\times N$ gives the number of control volumes per spatial direction, and corresponds to the number of degrees of freedom of the method. 
We conclude from the observed error rates that considering a high-order approximation will reach machine precision with a drastically reduced number of degrees of freedom. For example, we see that the $7^{th}$-order method is able to reach machine precision for a cell size as large as $\Delta x = L/32$. 

\subsection{Discontinuous magnetic field loop}
\begin{figure}
    \centering
    \includegraphics[width=\textwidth]{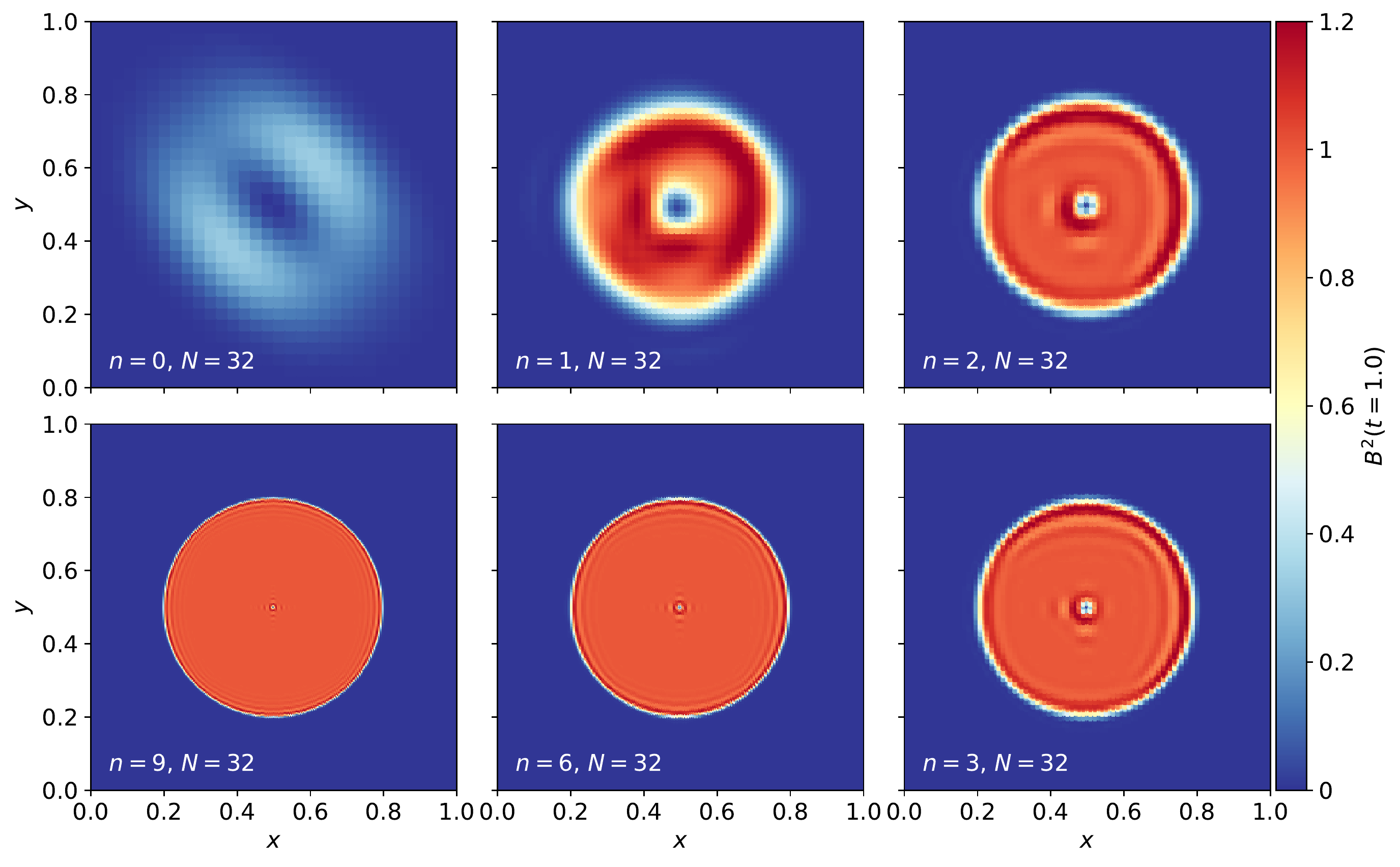}
    \caption{Discontinuous magnetic field loop advection test with increasing polynomial degree and $32$ cells on a side. \label{fig:sd-loop-order} }
\end{figure}

In this section we consider the initial conditions given by the discontinuous magnetic field loop test case, as introduced in section \ref{sec:overview}. We start by presenting in Fig.~\ref{fig:sd-loop-order} the solution maps computed at $t=1$ with the SD method while increasing the polynomial degree $n$, specifically for $n=0,1,2,3,6$ and $9$, and for $N=32$ cells per side. As we can see, increasing the order considerably improves the quality of the solution, and furthermore, even for a number of cells as small as $32$ per side, both the seventh- and tenth-order simulations ($n=6$ and $9$ respectively) show remarkable results preserving the shape of these discontinuous initial conditions. 

\begin{figure}
    \centering
    \includegraphics[width=\textwidth]{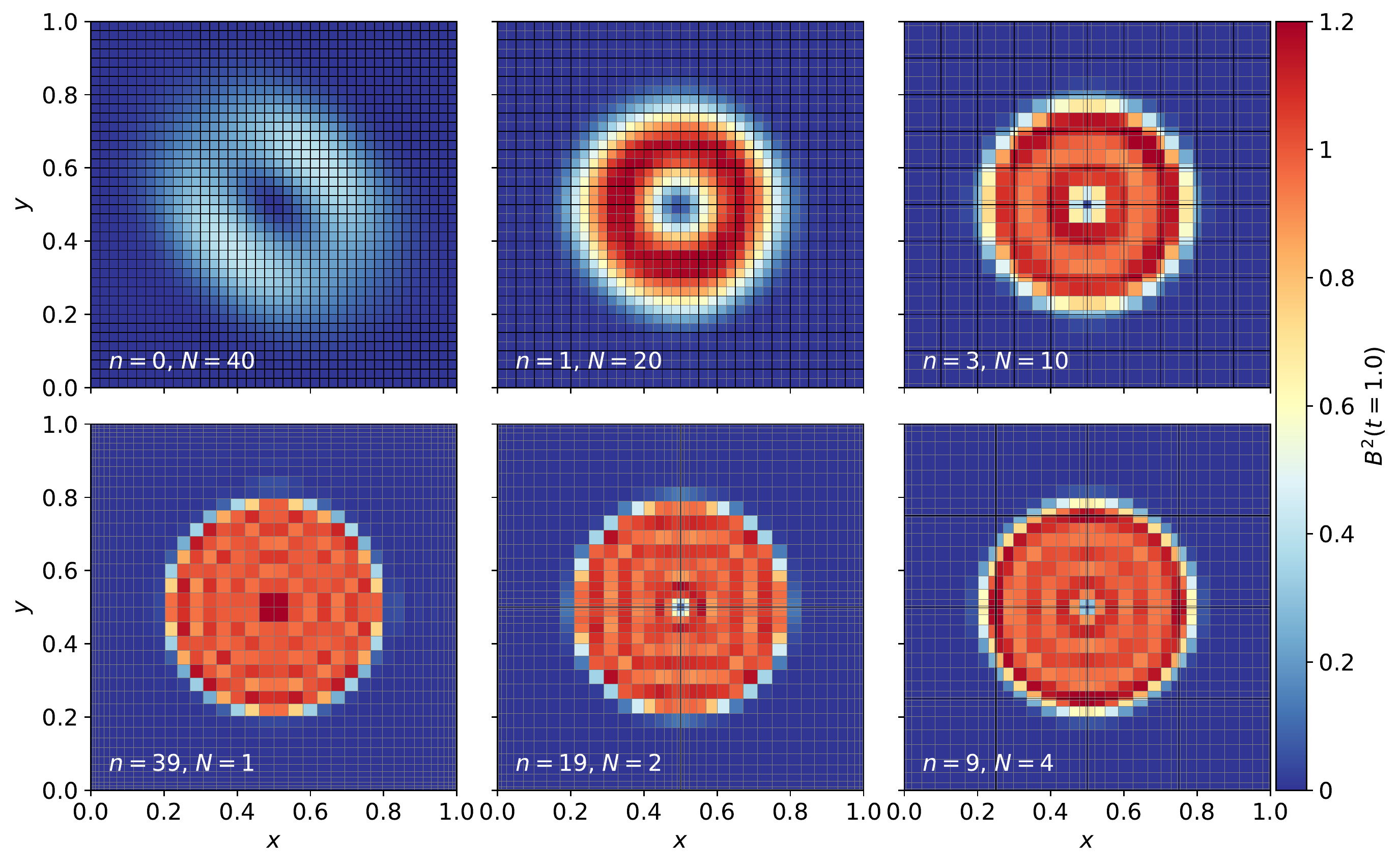}
    \caption{Discontinuous magnetic field loop advection test with increasing polynomial degree while maintaining the number of degrees of freedom equal to $40$. \label{fig:sd-loop-dof} }
\end{figure}

In Fig.~\ref{fig:sd-loop-dof} we go a step further, testing the "arbitrary high-order" character of our numerical implementation. In this figure we present again the solution maps at $t=1$, showing in black the mesh for the cells and in grey the mesh for the inner control volumes. While keeping constant the number of degrees of freedom per dimension $(n+1)\times N$, we show the increasingly better results as the order of the scheme is increased and the number of cells is decreased, keeping a constant product $(n+1)\times N=40$. In the most extreme case, of little practical interest, we go as far as testing a $40^{th}$-order method with one cell (as shown in the bottom-left panel of Fig.~\ref{fig:sd-loop-dof}). Surprisingly for us, this one cell simulation is able to preserve the initial conditions better than all the other cases. Indeed, in this extreme case, the flux points are "squeezed" towards the boundaries of the element, which results in an apparent loss of resolution of the control volumes at the centre of the element. The increased order of accuracy easily compensates for this effect.

\begin{figure}
    \centering
    \includegraphics[width=\textwidth]{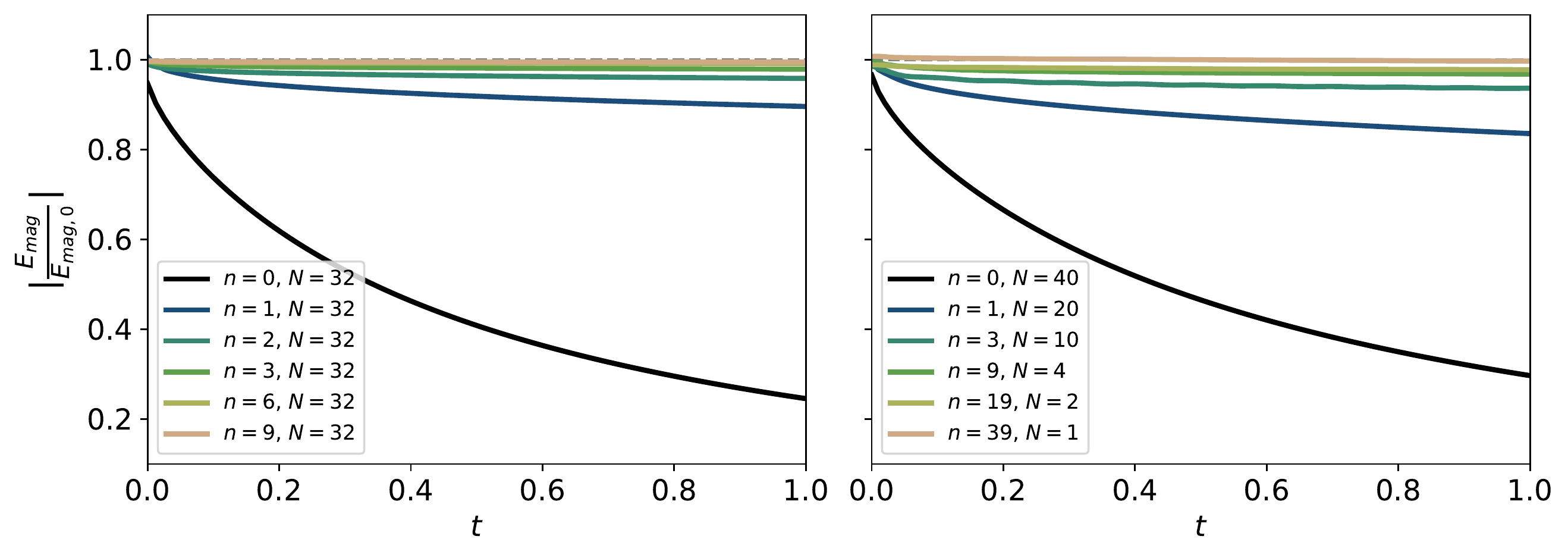}
    \caption{Normalised magnetic energy as a function of time for the discontinuous magnetic field loop advection test. The left panel shows the results for a fixed number of elements $N=32$ and an increasing order of accuracy, corresponding to Fig~\ref{fig:sd-loop-order}. The right panel shows the results for a fixed number of degrees of freedom $(n+1)\times N = 40$, corresponding to Fig~\ref{fig:sd-loop-dof}. \label{fig:sd-loop-EB} }
\end{figure}

We now present the performance of the method in preserving the magnetic energy. We show the normalised magnetic energy as a function of time for the simulations presented in Fig.~\ref{fig:sd-loop-order} (resp. Fig.~\ref{fig:sd-loop-dof}) on the left panel (resp. right panel) of Fig~\ref{fig:sd-loop-EB}. 
We see that going to higher order at fixed element resolution significantly improves the conservation property of the scheme. Our simulation with 32 elements and order 10 shows virtually no advection error anymore within the simulated time interval, 
at the expense of increasing the number of degrees of freedom  significantly. 
The second experiment, with a fixed number of degrees of freedom $(n+1)\times N = 40$, still shows 
a significant improvement in the energy conservation as the order of the method is increased. 
Our extreme case with only one element and a polynomial degree $n=39$ has also no visible advection errors in the magnetic energy evolution.
Note however that the computational cost of the method increases significantly with the order of accuracy, even when keeping the number of degrees of freedom constant. Regarding the conservation of the magnetic energy, for a given target accuracy, 
it is more efficient to go to higher order than to go to more computational elements.

\begin{figure}
   \includegraphics[width=0.94\textwidth,]{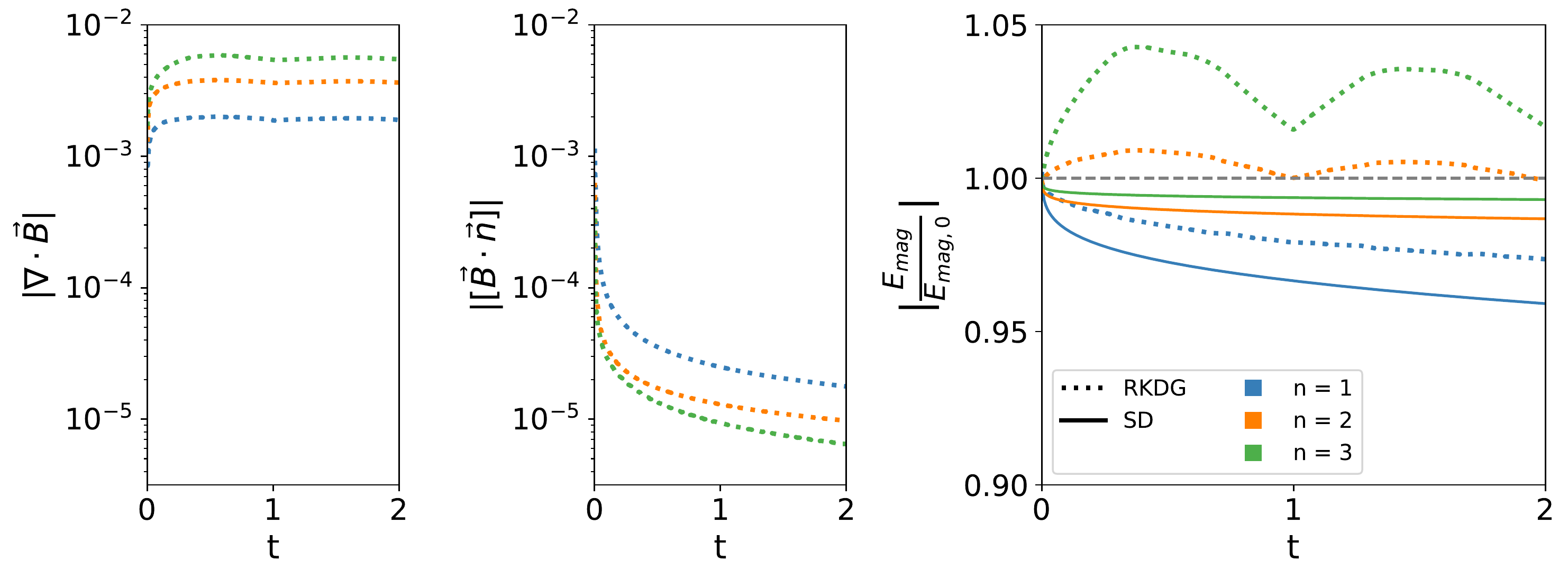}
   \begin{center}
   \includegraphics[width=1.0\textwidth]{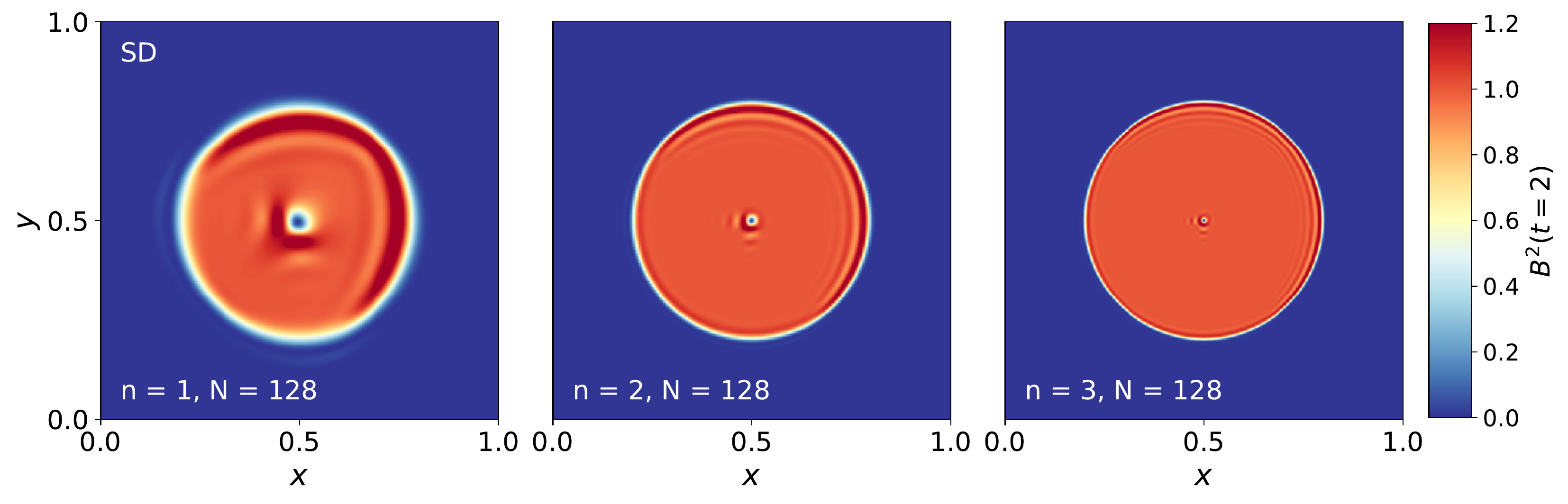}
    \caption{Same as Fig.~\ref{fig:rkdg-loopadvection-div-energy} but now for the SD scheme (solid lines). 
    For comparison, the results of the traditional RKDG scheme are shown as dotted lines. Note that the SD scheme has no local and no global divergence
    errors by construction (no solid lines in the top left and top middle panels).
    \label{fig:sd-loopadvection-div-energy}}
    \end{center}
\end{figure}

In order to compare our new SD scheme with the RKDG variants we have presented in section~\ref{sec:overview}, we show in Fig.~\ref{fig:sd-loopadvection-div-energy} the exact same field loop advection test for the SD implementation with $N=128$ elements per side. 
The reader is kindly asked to compare to Fig.~\ref{fig:rkdg-loopadvection-div-energy} through Fig.~\ref{fig:divc-loopadvection-div-energy}. 
The top left panel shows our results for the divergence errors of the numerical solution, compared to the traditional RKDG scheme, for both the volume and surface terms. 
This plot is meant as a joke, as obviously both terms are identically zero for the SD scheme, so only the traditional RKDG results are visible.
We confirm that the SD method preserves $\nabla\cdot\vec{B} = 0$ to machine precision, both in a global and in a local sense.
The right top panel shows again the magnetic energy evolution of the SD method, but this time with the same number of elements and order of accuracy than
the experiments performed in section~\ref{sec:overview}. We see that the SD method shows no spurious dynamo.
In the bottom panel of Fig.~\ref{fig:sd-loopadvection-div-energy}, we show the solution maps for magnetic energy density at $t=2$, 
in order to compare with the maps of Fig.~\ref{fig:rkdg-loopadvection-div-energy}, Fig.~\ref{fig:ldf-loopadvection-div-energy} and 
Fig. ~\ref{fig:divc-loopadvection-div-energy}. We note that the solution features a slight upwind asymmetry, 
as opposed to the solution of the DivClean RKDG method, especially for $n\leq3$. This upwind bias seems to disappear when moving to higher order.
A detailed comparison of the various schemes is presented in the next section.

\subsection{Rotating discontinuous magnetic field loop}
\label{subsection: rotating-hump}
In this section, we consider 
%the \textit{rotating hump} test case \cite{torrilhon2004},
the rotation of a discontinuous magnetic field loop. This test describes a linear velocity field $\vec{v} = (-y,x)^T$ acting on the magnetic field, resulting in a rotation around the origin. In this work, we use the following initial condition for the magnetic field $\vec{B_0}$:

\begin{equation}
\label{eq:magloop-ics}
    \vec{B}_0 = \begin{pmatrix} B_{x,0} \\ B_{y,0} \end{pmatrix} =  \begin{pmatrix} -A_0(y-y_c)/r \\ A_0(x-x_c)/r \end{pmatrix} \quad {\rm ~for~}r < r_0,
\end{equation} 
and $\vec{B}_0=0$ otherwise. We use here $A_0 = 0.001$, $r_0=\sfrac{1}{8}$ and $(x_c,y_c)=(\sfrac{3}{4}, \sfrac{1}{2})$. Then, the exact solution at time $t$ is given by:
\begin{equation}
    \vec{B}(\vec{x},t) = R(t)^{-1}B_0(R(t)\vec{x}),
\end{equation}
where $R(t)$ is a orthogonal matrix which rotates a vector by the angle $t$,
\begin{equation}
    R(t) = \begin{pmatrix} \cos(t) & -\sin(t) \\ \sin(t) & \cos(t) \end{pmatrix}.
\end{equation} 
Lastly, the computation domain considered is a box $[0,1]^2$ and at the boundary, the exact solution is prescribed in the ghost cells.

%This test case produces field lines circling around $\left( \sfrac{1}{2}, 0\right)$ and a distinct hump in $\norm{\vec{B}}$ with an essential radius of approximately $\sfrac{1}{8}$. In Figures \ref{fig:sd-rotatinghump-order} and \ref{fig:sd-rotatinghump-dof}, we show the solution computed by our proposed SD-ADER method varying the polynomial degree approximations and varying the polynomial degree approximation while keeping the total number of degrees of freedom constant, respectively. The solution is shown at $t=\pi$, corresponding to half a rotation. \todo{?}

\begin{figure}
    \centering
    \includegraphics[width=\textwidth]{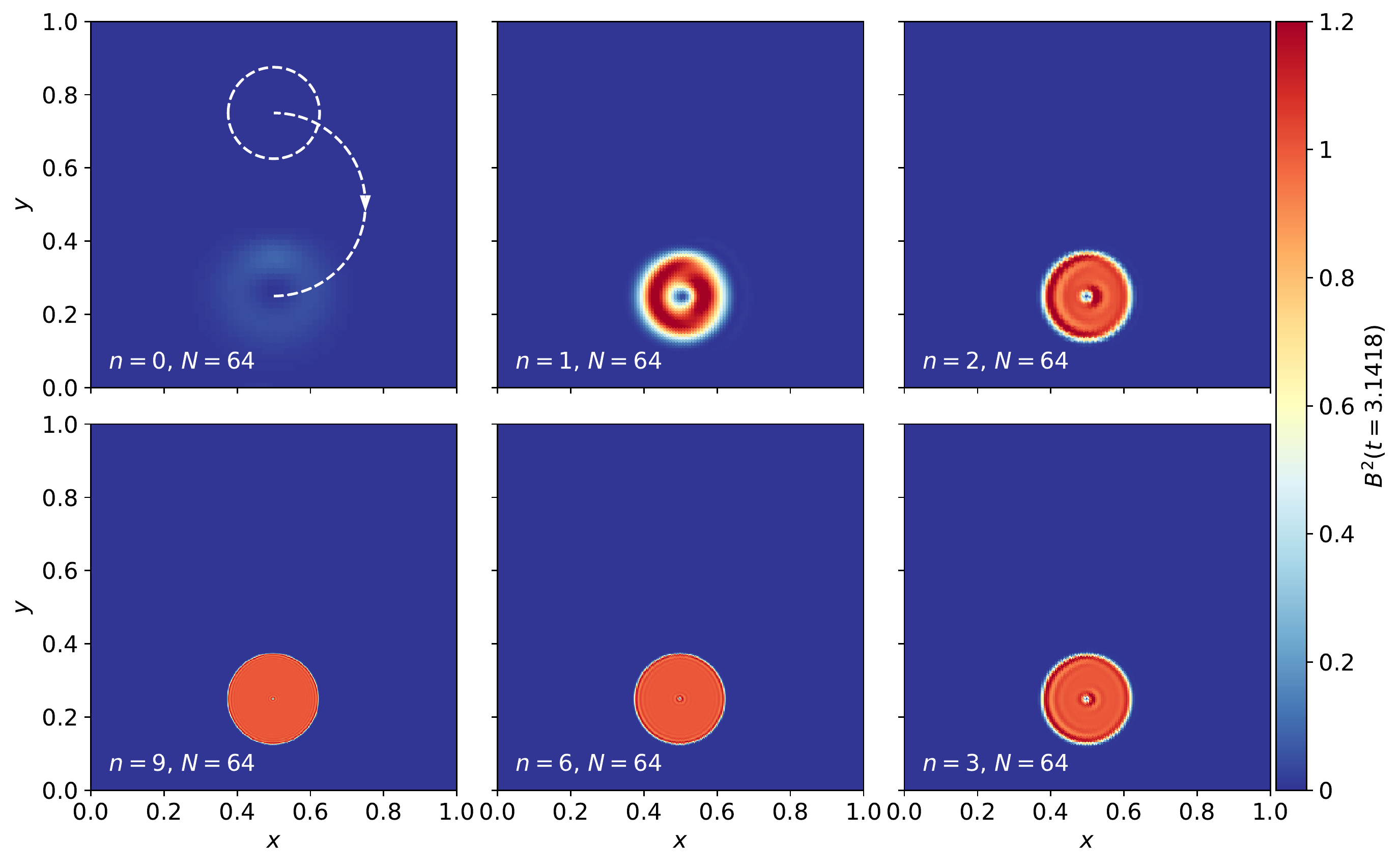}
    \caption{Rotating field loop test with increasing polynomial degree and $32$ cells on a side.
    %\todo{update figures}
    \label{fig:sd-rotatinghump-order} }
\end{figure}

\begin{figure}
    \centering
    \includegraphics[width=0.5\textwidth]{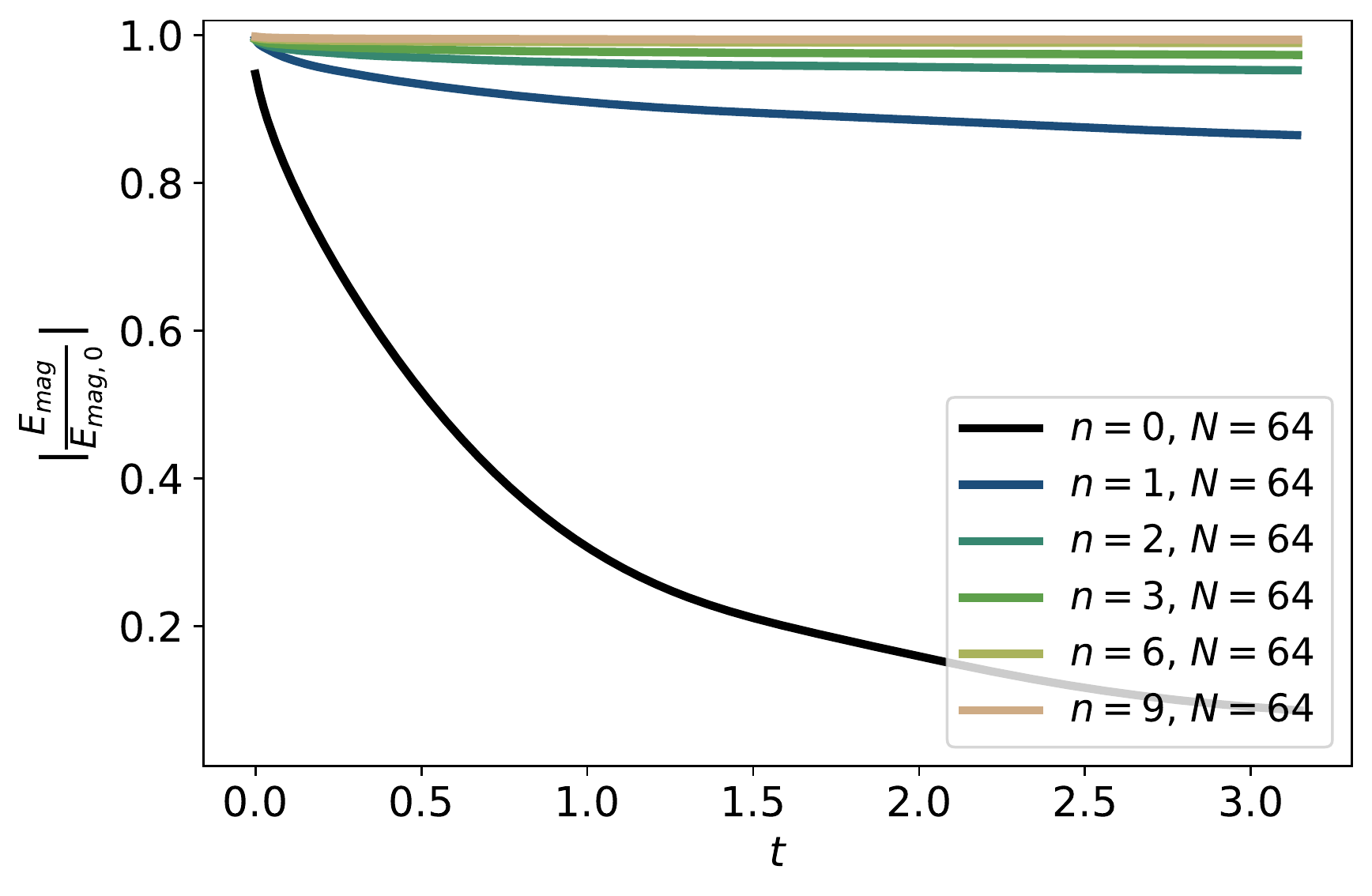}
    \caption{Normalised magnetic energy as a function of time for the rotating discontinuous field loop test. We show the results for a fixed number of elements $N=32$ and an increasing order of accuracy, corresponding to Fig~\ref{fig:sd-rotatinghump-order}.
    }
    \label{fig:sd-energy-rotation}
\end{figure}

In Fig.~ \ref{fig:sd-rotatinghump-order}, we show the solution computed by our proposed SD-ADER method varying the polynomial degree approximations. The solution is shown at $t=\pi$, corresponding to half a rotation. We observe that, as well as in the previous case, the method is able to preserve the discontinuous magnetic loop for $n$ greater or equal to 1. When comparing to Fig.~\ref{fig:sd-loop-order}, we have to highlight that the magnetic loop is being evolved up to a time $\pi$ times larger. Even then, the results remain similar, that is, increasingly better for higher order, thus showcasing the low numerical advection error that the method can reach. Furthermore, in Fig.~ \ref{fig:sd-energy-rotation}, the magnetic energy is shown. Once again, we expect the magnetic energy to remain constant over time, and indeed, we observe improvement in the conservation of the magnetic energy as the order is increased. In concrete, we observe for $n=6$ and $9$ a loss in magnetic energy below $1\%$. 

%\begin{figure}
%    \centering
%    \includegraphics[width=\textwidth]{chapters/img/mhd-loop-dof.pdf}
%    \caption{Rotating field loop test with increasing polynomial degree while maintaining the number of degrees of freedom equal to $40$.  
    %\todo{update figures} 
%    \label{fig:sd-rotatinghump-dof} }
%\end{figure}

%\begin{figure}
%    \centering
%    \includegraphics[width=\textwidth]{chapters/img/EB-dof.pdf}
%    \caption{Normalised magnetic energy as a function of time for the rotating field loop test. The left panel shows the results for a fixed number of elements $N=32$ and an increasing order of accuracy, corresponding to Fig~\ref{fig:sd-rotatinghump-order}. The right panel shows the results for a fixed number of degrees of freedom $(n+1)\times N = 40$, corresponding to Fig~\ref{fig:sd-rotatinghump-dof}. \label{fig:sd-rotatinghump-EB} 
    %\todo{update figures}
%    }
%\end{figure}

\section{Discussion}
\label{sec:mhd-discussion}

\subsection{Comparing SD to RKDG for the induction equation}

In this section, we compare in detail the different methods presented in this paper, namely our reference scheme, the traditional RKDG, a locally divergence-free basis variant of the scheme, called LDF, another variant of RKDG with divergence cleaning, called DivClean RKDG, and finally a novel Spectral Difference (SD) scheme specially designed for the induction equation, with the ADER time discretisation. The strong similarities with the Constrained Transport method would justify to call our new scheme using the long acronym CT-SD-ADER.

From a theoretical point of view, since the traditional RKDG scheme does not have any mechanism to deal with $\nabla\cdot\vec{B} \neq 0$, it is not so surprising to see this scheme perform relatively poorly. What is puzzling is why going to higher orders is so detrimental. Although the global contribution to the divergence error decreases with increasing order, the local divergence errors seem to increase with increasing order. As truncation errors decrease, the global divergence error decreases owing to smaller discontinuities at element boundaries, but the local divergence increases because of high-frequency and high-amplitude oscillations that damage the solution.
Considering a locally divergence-free polynomial basis for the magnetic field, as an explicit way to control the local divergence of the solution, seems like an obvious improvement of the scheme. However, we see that in this case the surface term, which measures the global divergence errors, becomes larger. We attribute this adverse effect to the fact that there are significantly less degrees of freedom available in the polynomial representation of the magnetic field, when comparing to the traditional RKDG scheme at the same order. 
Furthermore, as there is still no explicit mechanism to control global divergence errors, it is usually required to use the LDF basis in conjunction with an additional divergence cleaning mechanism to deal with the surface term. Indeed, we have shown that the divergence cleaning method (DivClean) provides an explicit, albeit non-exact, control on both the surface and the volume terms of the divergence errors, provided the two parameters, the hyperbolic cleaning speed $c_h$ and the diffusive coefficient $c_p^2$ are chosen appropriately.

With these considerations in mind, we designed a new numerical method based on the SD scheme, for which both the volume term and the surface term of the divergence errors vanish exactly. This new scheme satisfies an exact conservation law for the magnetic flux through a surface. We argue this is the natural way to interpret and solve the induction equation.
This approach, traditionally referred to as the Constrained Transport method, leads to a natural way to maintain zero divergence of $\vec{B}$ both locally and globally, as proved in Proposition \ref{proposition:pointwise_div_free} and Proposition \ref{proposition:globally_div_free}.

\begin{figure}
    \centering
    \includegraphics[width=0.94\textwidth]{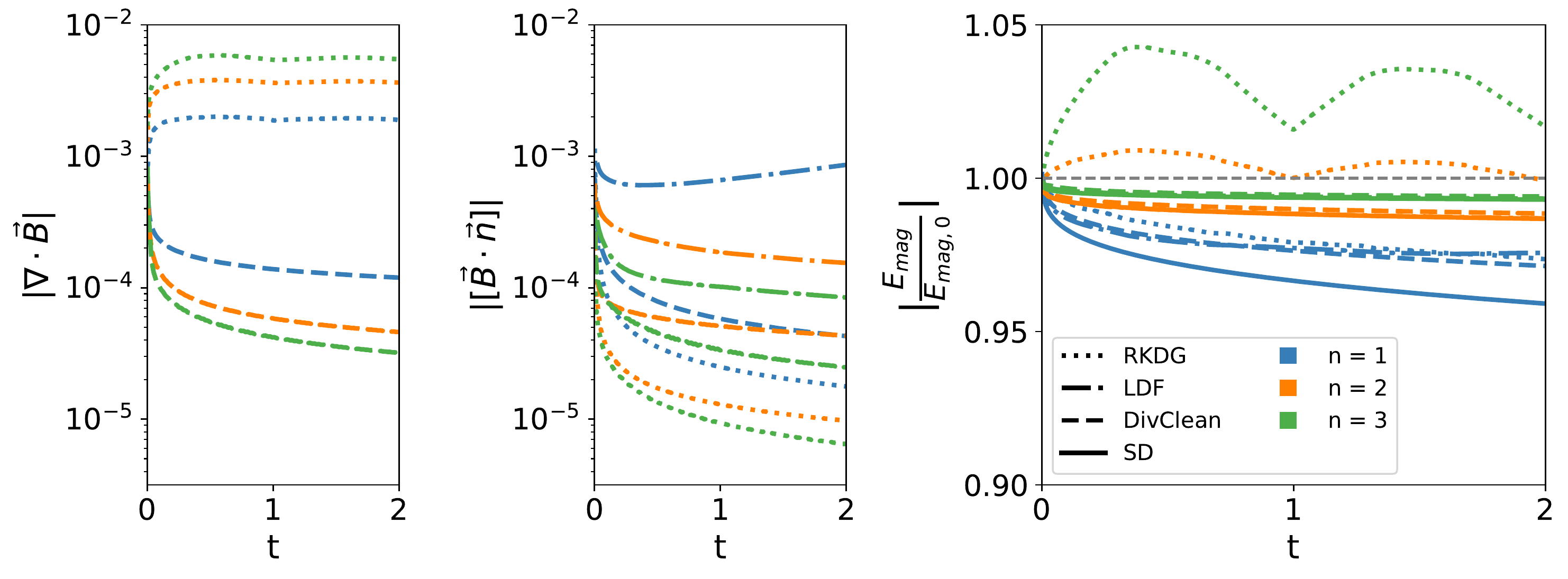}
    \caption{Local and global divergence errors and magnetic energy evolution of the four different methods discussed in this paper for the discontinuous magnetic field loop advection test
    and for polynomial degrees $n=1,~2,~3$. 
    \label{fig:all-div-energy} }
\end{figure}

We compared these 4 different methods by analyzing their performance when solving the advection of a discontinuous magnetic loop. 
The first (resp. second) panel of Fig.~\ref{fig:all-div-energy} shows the local (resp. global) divergence error of the schemes at different orders of accuracy. 
We note that for the SD scheme, we have zero contribution in both the volume and the surface terms. 
On the third panel of Fig.~\ref{fig:all-div-energy}, we show the magnetic energy evolution over time for the different methods.
The traditional RKDG method is the only one to exhibit a spurious dynamo at third and fourth orders.
The SD scheme appears slightly more diffusive than the other methods at second order, but its performance becomes comparable to LDF and DivClean at higher orders. 
Note that the extension to orders higher than $4$ for our new SD method is straightforward, as shown in the previous section, 
while the extension of the LDF method to orders higher than $4$ is quite cumbersome \citep[see for example][]{Guillet2019}.

\begin{figure}
    \centering
    \includegraphics[width=0.8\textwidth]{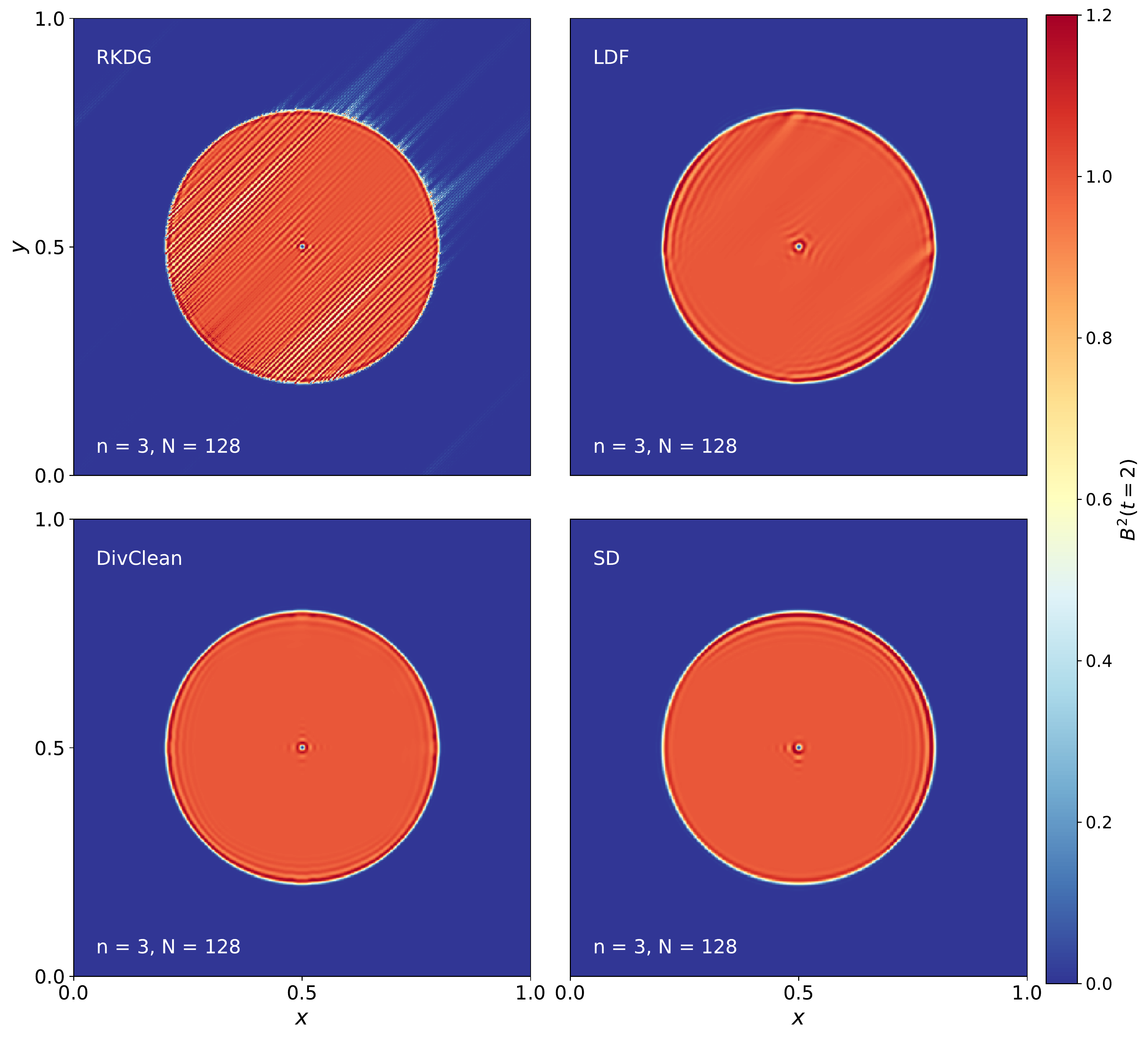}
    \caption{Maps of the magnetic energy density for the four different methods discussed in this paper for the discontinuous magnetic field loop advection test 
    and for polynomial degree $n=3$.
    \label{fig:disc-adv-allcomparison}
    }
\end{figure}

In Fig.~\ref{fig:disc-adv-allcomparison}, we show the maps of the magnetic energy density for the different schemes at fourth order and at $t=2$. 
First, we note that the magnetic energy distribution is well behaved for all the schemes, except RKDG, for which strong high-frequency oscillations are generated.
We also see that the solution computed using LDF retains some artifacts, which appear to be aligned with the velocity field. 
The solution computed with DivClean appears more symmetric and overall seems to have less artifacts, 
although some oscillations near the discontinuous boundary are still present, similarly to the solution computed with SD. 
To obtain the DivClean solution, some tuning of the parameters $c_h$ and $c_p$ is required. 
In particular, if $c_h$ is reduced from twice the advection velocity like here, to exactly equal to the advection velocity, 
the same artifacts that are seen in the solution computed with LDF appear in the solution using DivClean. 
It is also worth stressing again that the DivClean method comes with a price: a new equation and a new variable, whose physical interpretations are unclear.

A comparison of the methods with respect to their computational complexity is beyond the scope of this paper. In particular, the codes used to produce the numerical results have been developed with different programming languages and architectures. However, we can briefly comment on key similarities and differences between the DG-based methods presented and our SD-ADER method. We note that SD can be interpreted as a nodal, quadrature-free DG scheme \cite{May2011}, thus, making the proposed SD method not so different from a nodal DG one in terms of its computational complexity. Another key difference is the time-integration schemes used: for the DG-based schemes, we used SSP-RK time-integration whereas for the SD scheme we have used the ADER time-integration scheme. We note that to reach an $(n+1)$-order approximation in time, the ADER algorithm requires $n+1$ flux evaluations per time slice \cite{Jackson2017}, yielding an overall complexity of $(n+1)^2$ in time. Then, it becomes computationally more expensive than an explicit RK scheme, as the number of stages needed to reach an $(n+1)$-order approximation is typically well below $(n+1)^2$. However, as noted in \cite{Dumbser2018}, the ADER procedure can be formulated as a completely local predictor step suited for vectorisation, reducing then the complexity to $n+1$, whereas the RK scheme requires communication with its neighbours at every stage.

\subsection{SD method for non-trivial velocity fields}
In subsection \ref{subsection: rotating-hump}, we consider the problem of a rotating velocity field. We show the ability of our method to solve problems with non-trivial velocity fields, as well as
Dirichlet boundary conditions. For approximation polynomial degree of $n=1$, we obtain similar qualitative results to those of \cite{Torrilhon2004} (given that the initial $B_0$ is different). As we increase the approximation order, we can observe that the numerical solution converges to the analytical one.

\subsection{Extension of our new method to three space dimensions}

In this section, we speculate about a possible straightforward extension of our scheme to three dimensions in space. It is not the scope of this paper to present a detailed implementation of this algorithm, however, we want to stress that this extension is not only possible, but also relatively easy and consistent with the present work. They are however a few key differences with respect to the 2D case. 
The first difference comes from the definition of the magnetic flux and from the resulting magnetic field. 
We now define the magnetic flux of $B_x$ across a rectangle sitting in the plane $x=x_i$ and 
defined by the 4 points $(0,0)$, $(0,z_k)$, $(y_j,0)$ and $(y_j,z_k)$ as:
\begin{equation}
\phi_x(x_i,y_j,z_k) = \int_0^{y_j} \int_0^{z_k} B_x(x_i,y,z) {\rm d}y {\rm d}z,
\end{equation} 
where the coordinates $y_j$ and $z_k$ correspond to the flux points in each direction, or in other words, to the corner points of each control volume inside the element.

The magnetic flux is then interpolated everywhere inside the element using Lagrange polynomials defined using the flux points.
\[\phi_x(x,y,z) = \sum_{i=0}^{n+1} \sum_{j=0}^{n+1} \sum_{k=0}^{n+1}\phi_x(x_i,y_j,z_k) \ell_i(x) \ell_j(y) \ell_k(z). \] 
The magnetic field inside the element is obtained through a second-order derivative as follows:
\begin{equation}
B_{x}(x,y,z) = \partial^2_{yz} \phi_x = \sum_{i=0}^{n+1} \sum_{j=0}^{n+1} \sum_{k=0}^{n+1} 
\phi_x(x_i,y_j,z_k) \ell_i(x) \ell^{\prime}_j(y) \ell^{\prime}_k(z).
\end{equation}
Note that this formulation is equivalent to the alternative approach we describe below using the vector potential.
It is however important to understand that this interpolation is used only at initialisation, to make sure the corresponding
magnetic field is strictly divergence free.

The next step is to evaluate the magnetic field at the solution points, which, in the 3D case, are now located at the centre of the face in the staggered direction: 
$B_x(x^f_i,y^s_j,z^s_k)$, $B_y(x^s_i,y^f_j,z^s_k)$ and $B_z(x^s_i,y^s_j,z^f_k)$, 
where the $f$ and $s$ superscripts correspond again to flux and solution points respectively. 
Once the field has been initialised on the solution points, we then interpolate the field within each face of the control volumes using Lagrange polynomials, which are defined using the solution points as in the traditional SD method.
Using these definitions, it is straightforward to generalise Proposition~\ref{proposition:pointwise_div_free} to the 3D case, and prove that $\nabla\cdot\vec{B}=0$.

The components of the electric field are defined at the centre of the edges between control volumes, located at flux points in the directions orthogonal to the component, and at solution points along the components direction: $E_x(x^s_i,y^f_j,z^f_k)$, $E_y(x^f_i,y^s_j,z^f_k)$ and $E_z(x^f_i,y^f_j,z^s_k)$. The electric field is again defined as $\vec{E}= - \vec{v}\times\vec{B}$, therefore this method requires to know the orthogonal velocities at those same edges,
and to solve a 1D Riemann problem at element's faces and a 2D Riemann problem at element's edges.
As in Proposition~\ref{proposition:globally_div_free}, the SD update of the magnetic field is obtained directly using a pointwise update at the magnetic field solution points:
\begin{equation}
\partial_t B_{x} = \partial_z E_y - \partial_y E_z,~~~ 
\partial_t B_{y} = \partial_x E_z - \partial_z E_x~~~{\rm and}~~~ 
\partial_t B_{z} = \partial_y E_x - \partial_x E_y.
\end{equation}

It follows trivially, like in the 2D case, that
\begin{equation}
\partial_t \left( \partial_x B_{x} + \partial_y B_{y} + \partial_z B_{z}\right) = 0.
\end{equation}

We have here again an equivalence between the SD method applied to the magnetic field and a similar SD method 
applied to the vector potential. It is however more difficult in 3D to compute the vector potential from the magnetic field. 
It requires a complex inversion and the choice of a gauge, using for example the Coulomb gauge, for which $\nabla \cdot \vec{A}=0$.

Assuming we know the vector potential, we define for each component the line integral over the component's direction, as shown here for the $z$-direction:
\begin{equation}
\alpha_z(x_i,y_j,z_k) =\int_0^{z_k} A_z(x_i,y_j,z) {\rm d}z.
\end{equation} 

As for the magnetic flux, this quantity is defined at the corner points of the control volumes using flux points in each direction. 
We can then use the Lagrange polynomials defined using the flux points to compute the vector potential everywhere as:
\begin{equation}
A_{z}(x,y,z) = \partial_z \alpha_z = \sum_{i=0}^{n+1} \sum_{j=0}^{n+1} \sum_{k=0}^{n+1} \alpha_z(x_i,y_j,z_k) \ell_i(x) \ell_j(y) \ell^{\prime}_k(z).
\end{equation}

We can now evaluate the vector potential at the corresponding solution points, which are, as for $\vec{E}$, 
defined at the edges of the control volumes:  $A_x(x^s_i,y^f_j,z^f_k)$, $A_y(x^f_i,y^s_j,z^f_k)$ and $A_z(x^f_i,y^f_j,z^s_k)$. 

Once we know the polynomial representation of the vector potential, the magnetic field can be derived using pointwise 
derivatives and $\vec{B} = \nabla \times \vec{A}$. 
The vector potential can finally be updated directly at its solution points using (shown here only for $A_z$):
\begin{equation}
\partial_t A_z = -v_x \partial_x A_z - v_y \partial_y A_z + v_x \partial_z A_x + v_y \partial_z A_y. 
\end{equation}

This is again the vector potential equation, although in a more complex form than in the 2D case. 
It can however be solved using our SD scheme, exactly like in 2D.

%\revA{I am a bit confused about why we only introduce the vector potential much later on when we talk about 3d, not in the beginning like we do in 2-d.}

\subsection{Extension of our new method to ideal MHD }

The natural progression of this work is to extend the proposed SD method to the full magneto-hydrodynamics equations. The first difficulty is to solve 2D Riemann problems at element edges. Fortunately, 2D Riemann solvers in the context of ideal MHD have been already developed in the past years in multiple implementations of Constrained  Transport for the FV Godunov method \cite{Londrillo2004,Teyssier2007,balsara2010,balsara2012,Balsara2014,balsara2015a,balsara2017}.

As for the time stepping, the ADER methodology is trivially extended to 3-D and nonlinear problems \cite[see e.g.][]{dumbser_ader_2013}. Our proposed version of ADER only differs in the fact that we do not remain local during the iterative process, as we require Riemann solvers as part of the SD space discretization. This means that in the MHD case, we ought to use an appropriate Riemann solver as described above.

The second difficulty comes from finding the appropriate shock capturing techniques for the SD method, which traditionally has been achieved through artificial viscosity \cite{Premasuthan2014}. Finding both a way to enforce preservation of positiveness and not clipping smooth extrema, while constraining as least as possible the performance of the method, is of extreme importance. Recent advances in shock capturing methods, such as \cite{Vilar2019}, provide a natural way of performing sub-cell limiting in a nodal discontinuous Galerkin method based on the \textit{a posteriori} subcell limiting strategy (MOOD) \cite{Dumbser2016} that can guarantee both positivity and as little as possible extrema clipping in smooth profiles. This methodology seems promising to be applied to our SD method in the context of the ideal MHD equations when used in combination with a robust finite volume scheme that preserves the divergence free nature of the solution.
\section{Conclusions}
\label{sec:conclusion}

In this work, we have analysed in detail several variants of the high-order DG method with RK time integration for the induction 
equation, while attempting to preserve the constraint $\nabla\cdot\vec{B}=0$ with various degrees of success. 
We have then presented a novel, arbitrary high-order numerical scheme based on a modification of the Spectral Difference (SD) method with ADER time integration for the induction equation. This new scheme preserves $\nabla\cdot\vec{B}=0$ exactly by construction. 
It is a natural extension of the Constrained Transport scheme to the SD method. We have proved that both the volume term and the surface term in the norm definition of the divergence vanish. 
We have also reformulated our scheme in terms of the vector potential, which allows a direct connection with a classical 
SD method for the evolution of the vector potential, allowing us to analyse its stability and dispersion properties, with results similar to \cite{abeele2008}. Furthermore, we show that the combination of ADER and SD result in a stable method when choosing the appropriate CFL condition.

We have shown with various numerical experiments that our method converges at the expected order, namely $\Delta x^{n+1}$, 
where $n$ is the polynomial degree of the adopted interpolation Lagrange polynomials and $\Delta x$ the element size. 
We have also considered the discontinuous field loop advection test case \cite{Gardiner2005}, a problem known to reveal artifacts caused by not preserving $\nabla\cdot\vec{B}=0$. We have shown again that our new method behaves well, up to incredibly high orders (polynomial degree $n=39$), conserving the magnetic energy almost exactly by drastically reducing advection errors, provided the order is high enough. Furthermore, we also test our method using a non-trivial velocity field and show our method leads to the correct solution, and that we qualitatively get similar results as in \cite{Torrilhon2004}.

We have then compared our novel method with the high-order DG variants presented in the first part of the paper. The magnetic energy evolution and the solution maps of the SD-ADER scheme all show qualitatively similar and overall good performances when compared to the Divergence Cleaning method applied to RKDG, but without the need for an additional equation and an extra variable 
to help controlling the divergence errors. We have finally discussed our future plans to extend this work to three dimensions and to fully non-linear ideal MHD.
\section{Acknowledgments}
We gratefully thank R. Abgrall (University of Zurich) and S. Mishra (ETH Zurich) for the fruitful discussions and insights regarding this work. This research was supported in part through computational resources provided by ARC-TS (University of Michigan) and CSCS, the Swiss National Supercomputing Centre. %MHV acknowledges financial support from MIDAS. DAVR acknowledges funding from/is supported by ?GRANT?.
\bibliographystyle{unsrt}  
\bibliography{references} 

\newpage
\section{Appendix}
\subsection{Locally divergence-free basis}
\label{ap:LDF}
In order to design a locally-divergence free basis to represent $\vec{B}$ up to order $n+1$, the vector basis elements are computed as the curl of the elements of the polynomial space $\mathbb{P}^{n+1}$, which contains all polynomials in $x$ and $y$ up to degree $n+1$. Furthermore, as noted in \cite{klingenberg2017}, an orthogonal basis yields better conditioned mass matrices, so we apply the Gram-Schmidt orthogonalisation algorithm with the inner product:
\[\eta(\vec{b}_i, \vec{b}_j) = \int_{[-1,1]^2} \vec{b}_i \cdot \vec{b}_j {\rm d} x {\rm d}y.\]
The orthogonal and normalised basis vectors (up to $4^{th}$ order of approximation) are given below. These were obtained through the symbolic computation package {\ttfamily{sympy}} \cite{sympy}.
\begin{gather*}
 \mhdbasis^1 = {\rm span}\left(\left\{ 
 \begin{pmatrix} 1.0 \\ 0 \end{pmatrix}, 
 \begin{pmatrix} 0 \\ 1.0 \end{pmatrix},
 \begin{pmatrix} \sqrt{3}y \\ 0  \end{pmatrix},
 \begin{pmatrix} 0\\  \sqrt{3}x  \end{pmatrix},
 \begin{pmatrix} \sqrt{\frac{3}{2}}x \\ -\sqrt{\frac{3}{2}}y \end{pmatrix}\right\}\right)
\end{gather*}

\begin{equation*}
\begin{split}
 \mhdbasis^2 = \mhdbasis^1 \cup {\rm span}\left(\Bigg\{
 \sqrt{30}\begin{pmatrix} \frac{3x^2-1}{12} \\ -\frac{xy}{2} \end{pmatrix},
 \sqrt{30}\begin{pmatrix} -\frac{xy}{2} \\ \frac{3y^2-1}{12} \end{pmatrix}, \sqrt{5} \begin{pmatrix} \frac{3y^2 - 1}{2} \\ 0\end{pmatrix}, \right. \left. \sqrt{5} \begin{pmatrix} 0 \\ \frac{3x^2-1}{2}\end{pmatrix}
 \Bigg\}\right)
 \end{split}
\end{equation*}

\begin{equation*}
\begin{split}
\mhdbasis^3 = \mhdbasis^2 \cup {\rm span}\left(\Bigg\{
\frac{\sqrt{42}\sqrt{83}}{166}\begin{pmatrix} 5x^3-4x \\ -15yx^2+4y \end{pmatrix}, \frac{\sqrt{30}}{4}\begin{pmatrix}3x^2y-y \\ -3y^2x + x \end{pmatrix}, \frac{\sqrt{7}}{2}\begin{pmatrix} 5y^3-3y \\ 0 \end{pmatrix},
\frac{\sqrt{7}}{2}\begin{pmatrix}0\\ 5x^3-3x \end{pmatrix}, \right.\\ \left.
\frac{\sqrt{165585}}{1824}\begin{pmatrix} -\frac{56x^3}{83}-2x(12y^2-1)+\frac{410x}{83} \\ 8y^3 + \frac{14y(12x^2-1.0)}{83} - \frac{562y}{83} \end{pmatrix}\Bigg\}\right).
\end{split}
\end{equation*}
A detailed discussion on the approximation properties of this polynomial vector space can be found in \cite{Cockburn2004}.

\subsection{Divergence cleaning}
\label{ap:DivClean}
%\todo[inline]{Maybe this goes to appendix or not at all.}
We evolve the system defined by Eq.~\eqref{eq:glm-induction-eq} as described in \cite{klingenberg2017} and using the following steps:
\begin{enumerate}
    \item We apply SSP-RK to the DG discretisation of the induction equation in its divergence form as in Eq.~\eqref{eq:mhd-induction-eq-div}.
    \item{ We then apply in an operator split fashion SSP-RK to the DG discretisation of the system
    \begin{equation*}
    \begin{split}
    \partial_t \vec{B} + \nabla \psi &= 0,\\
    \partial_t \psi + c_h^2\nabla\cdot\vec{B} &= 0.
    \end{split}
    \end{equation*}
    }
    \item{We finally apply operator splitting to the source term of the parabolic term
    \[ \psi^{n+1} :=\exp\left(-\frac{c_h^2}{c_p^2}\Delta t\right) \psi^{n+1/2}. \]
    }
\end{enumerate}

\end{document}